
\newcount\mgnf\newcount\tipi\newcount\tipoformule\newcount\greco

\tipi=2          
\tipoformule=0   


\global\newcount\numsec
\global\newcount\numfor
\global\newcount\numtheo
\global\advance\numtheo by 1

\def\senondefinito#1{\expandafter\ifx\csname#1\endcsname\relax}

\def\SIA #1,#2,#3 {\senondefinito{#1#2}%
\expandafter\xdef\csname #1#2\endcsname{#3}\else
\write16{???? ma #1,#2 e' gia' stato definito !!!!} \fi}

\def\etichetta(#1){(\veroparagrafo.\veraformula)%
\SIA e,#1,(\veroparagrafo.\veraformula) %
\global\advance\numfor by 1%
\write15{\string\FU (#1){\equ(#1)}}%
\write16{ EQ #1 ==> \equ(#1) }}

\def\letichetta(#1){\veroparagrafo.\verotheo
\SIA e,#1,{\veroparagrafo.\verotheo}
\global\advance\numtheo by 1
\write15{\string\FU (#1){\equ(#1)}}
\write16{ Sta \equ(#1) == #1 }}

\def\tetichetta(#1){\veroparagrafo.\veraformula 
\SIA e,#1,{(\veroparagrafo.\veraformula)}
\global\advance\numfor by 1
\write15{\string\FU (#1){\equ(#1)}}
\write16{ tag #1 ==> \equ(#1)}}

\def\FU(#1)#2{\SIA fu,#1,#2 }

\def\etichettaa(#1){(A\veroparagrafo.\veraformula)%
\SIA e,#1,(A\veroparagrafo.\veraformula) %
\global\advance\numfor by 1%
\write15{\string\FU (#1){\equ(#1)}}%
\write16{ EQ #1 ==> \equ(#1) }}

\def\BOZZA{
\def\alato(##1){%
 {\rlap{\kern-\hsize\kern-1.4truecm{$\scriptstyle##1$}}}}%
\def\aolado(##1){%
 {
{
 \rlap{\kern-1.4truecm{$\scriptstyle##1$}}}}}
 }

\def\alato(#1){}
\def\aolado(#1){}

\def\veroparagrafo{\number\numsec}
\def\veraformula{\number\numfor}
\def\verotheo{\number\numtheo}

\def\Eq(#1){\eqno{\etichetta(#1)\alato(#1)}}
\def\eq(#1){\etichetta(#1)\alato(#1)}
\def\leq(#1){\leqno{\aolado(#1)\etichetta(#1)}}
\def\teq(#1){\tag{\aolado(#1)\tetichetta(#1)\alato(#1)}}
\def\Eqa(#1){\eqno{\etichettaa(#1)\alato(#1)}}
\def\eqa(#1){\etichettaa(#1)\alato(#1)}
\def\eqv(#1){\senondefinito{fu#1}$\clubsuit$#1
\write16{#1 non e' (ancora) definito}%
\else\csname fu#1\endcsname\fi}
\def\equ(#1){\senondefinito{e#1}\eqv(#1)\else\csname e#1\endcsname\fi}

\def\Lemma(#1){\aolado(#1)Lemma \letichetta(#1)}%
\def\Theorem(#1){{\aolado(#1)Theorem \letichetta(#1)}}%
\def\Proposition(#1){\aolado(#1){Proposition \letichetta(#1)}}%
\def\Corollary(#1){{\aolado(#1)Corollary \letichetta(#1)}}%
\def\Remark(#1){{\noindent\aolado(#1){\bf Remark \letichetta(#1).}}}%
\def\Definition(#1){{\noindent\aolado(#1){\bf Definition 
\letichetta(#1)$\!\!$\hskip-1.6truemm}}}
\def\Example(#1){\aolado(#1) Example \letichetta(#1)$\!\!$\hskip-1.6truemm}

\def\include#1{
\openin13=#1.aux \ifeof13 \relax \else
\input #1.aux \closein13 \fi}

\openin14=\jobname.aux \ifeof14 \relax \else
\input \jobname.aux \closein14 \fi
\openout15=\jobname.aux


{\count255=\time\divide\count255 by 60 \xdef\hourmin{\number\count255}
        \multiply\count255 by-60\advance\count255 by\time
   \xdef\hourmin{\hourmin:\ifnum\count255<10 0\fi\the\count255}}

\def\oramin{\hourmin }

\def\data{\number\day/\ifcase\month\or january \or february \or march \or april
\or may \or june \or july \or august \or september
\or october \or november \or december \fi/\number\year;\ \oramin}

\newcount\pgn \pgn=1
\def\foglio{\number\numsec:\number\pgn
\global\advance\pgn by 1}
\def\foglioa{A\number\numsec:\number\pgn
\global\advance\pgn by 1}

\footline={\rlap{\hbox{\copy200}}\hss\tenrm\folio\hss}

\def\TIPIO{
\font\setterm=amr7 
\def \settepunti{\def\rm{\fam0\setterm}
\textfont0=\setterm   
\normalbaselineskip=9pt\normalbaselines\rm }\let\nota=\settepunti}

\def\TIPITOT{
\font\twelverm=cmr12
\font\twelvei=cmmi12
\font\twelvesy=cmsy10 scaled\magstep1
\font\twelveex=cmex10 scaled\magstep1
\font\twelveit=cmti12
\font\twelvett=cmtt12
\font\twelvebf=cmbx12
\font\twelvesl=cmsl12
\font\ninerm=cmr9
\font\ninesy=cmsy9
\font\eightrm=cmr8
\font\eighti=cmmi8
\font\eightsy=cmsy8
\font\eightbf=cmbx8
\font\eighttt=cmtt8
\font\eightsl=cmsl8
\font\eightit=cmti8
\font\sixrm=cmr6
\font\sixbf=cmbx6
\font\sixi=cmmi6
\font\sixsy=cmsy6
\font\twelvetruecmr=cmr10 scaled\magstep1
\font\twelvetruecmsy=cmsy10 scaled\magstep1
\font\tentruecmr=cmr10
\font\tentruecmsy=cmsy10
\font\eighttruecmr=cmr8
\font\eighttruecmsy=cmsy8
\font\seventruecmr=cmr7
\font\seventruecmsy=cmsy7
\font\sixtruecmr=cmr6
\font\sixtruecmsy=cmsy6
\font\fivetruecmr=cmr5
\font\fivetruecmsy=cmsy5
\textfont\truecmr=\tentruecmr
\scriptfont\truecmr=\seventruecmr
\scriptscriptfont\truecmr=\fivetruecmr
\textfont\truecmsy=\tentruecmsy
\scriptfont\truecmsy=\seventruecmsy
\scriptscriptfont\truecmr=\fivetruecmr
\scriptscriptfont\truecmsy=\fivetruecmsy
\def \eightpoint{\def\rm{\fam0\eightrm}
\textfont0=\eightrm \scriptfont0=\sixrm \scriptscriptfont0=\fiverm
\textfont1=\eighti \scriptfont1=\sixi   \scriptscriptfont1=\fivei
\textfont2=\eightsy \scriptfont2=\sixsy   \scriptscriptfont2=\fivesy
\textfont3=\tenex \scriptfont3=\tenex   \scriptscriptfont3=\tenex
\textfont\itfam=\eightit  \def\it{\fam\itfam\eightit}%
\textfont\slfam=\eightsl  \def\sl{\fam\slfam\eightsl}%
\textfont\ttfam=\eighttt  \def\tt{\fam\ttfam\eighttt}%
\textfont\bffam=\eightbf  \scriptfont\bffam=\sixbf
\scriptscriptfont\bffam=\fivebf  \def\bf{\fam\bffam\eightbf}%
\tt \ttglue=.5em plus.25em minus.15em
\setbox\strutbox=\hbox{\vrule height7pt depth2pt width0pt}%
\normalbaselineskip=9pt
\let\sc=\sixrm  \let\big=\eightbig  \normalbaselines\rm
\textfont\truecmr=\eighttruecmr
\scriptfont\truecmr=\sixtruecmr
\scriptscriptfont\truecmr=\fivetruecmr
\textfont\truecmsy=\eighttruecmsy
\scriptfont\truecmsy=\sixtruecmsy }\let\nota=\eightpoint}

\newfam\msbfam   
\newfam\truecmr  
\newfam\truecmsy 
\newskip\ttglue
\ifnum\tipi=0\TIPIO \else\ifnum\tipi=1 \TIPI\else \TIPITOT\fi\fi

\def\e{\epsilon}

\def\E{{I\kern-.25em{E}}}
\def\N{{I\kern-.25em{N}}}
\def\M{{I\kern-.25em{M}}}
\def\R{{I\kern-.25em{R}}}
\def\Z{{Z\kern-.425em{Z}}}
\def\1{{1\kern-.25em\hbox{\rm I}}}
\def\eu{{1\kern-.25em\hbox{\sm I}}}

\def\C{{I\kern-.64em{C}}}
\def\P{{I\kern-.25em{P}}}
\def\eop{{ \vrule height7pt width7pt depth0pt}\par\bigskip}



\def\chap #1#2{\line{\ch #1\hfill}\numsec=#2\numfor=1}

\def\ba{{\backslash}}

\def\sqr#1#2{{\vcenter{\vbox{\hrule height.#2pt
     \hbox{\vrule width.#2pt height#1pt \kern#1pt
   \vrule width.#2pt}\hrule height.#2pt}}}}
\def\qed{ $\mathchoice\sqr64\sqr64\sqr{2.1}3\sqr{1.5}3$} 


\newcount\foot
\foot=1
\def\note#1{\footnote{${}^{\number\foot}$}{\ftn #1}\advance\foot by 1}
\def\tag #1{\eqno{\hbox{\rm(#1)}}}
\def\frac#1#2{{#1\over #2}}

\def\text#1{\quad{\hbox{#1}}\quad}
\def\newpage{\vfill\eject}

\def\thanks{\noindent{\bf Aknowledgements: }}



\font\ch=cmbx12

\font\ftn=cmr8

\font\it=cmti10
\font\sm=cmr7

%
\catcode`\X=12\catcode`\@=11
\def\n@wcount{\alloc@0\count\countdef\insc@unt}
\def\n@wwrite{\alloc@7\write\chardef\sixt@@n}
\def\n@wread{\alloc@6\read\chardef\sixt@@n}
\def\crossrefs#1{\ifx\alltgs#1\let\tr@ce=\alltgs\else\def\tr@ce{#1,}\fi
   \n@wwrite\cit@tionsout\openout\cit@tionsout=\jobname.cit 
   \write\cit@tionsout{\tr@ce}\expandafter\setfl@gs\tr@ce,}
\def\setfl@gs#1,{\def\@{#1}\ifx\@\empty\let\next=\relax
   \else\let\next=\setfl@gs\expandafter\xdef
   \csname#1tr@cetrue\endcsname{}\fi\next}
\newcount\sectno\sectno=0\newcount\subsectno\subsectno=0\def\r@s@t{\relax}
\def\resetall{\global\advance\sectno by 1\subsectno=0
  \gdef\firstpart{\number\sectno}\r@s@t}
\def\resetsub{\global\advance\subsectno by 1
   \gdef\firstpart{\number\sectno.\number\subsectno}\r@s@t}
\def\v@idline{\par}\def\firstpart{\number\sectno}
\def\l@c@l#1X{\firstpart.#1}\def\gl@b@l#1X{#1}\def\t@d@l#1X{{}}
\def\m@ketag#1#2{\expandafter\n@wcount\csname#2tagno\endcsname
     \csname#2tagno\endcsname=0\let\tail=\alltgs\xdef\alltgs{\tail#2,}%
  \ifx#1\l@c@l\let\tail=\r@s@t\xdef\r@s@t{\csname#2tagno\endcsname=0\tail}\fi
   \expandafter\gdef\csname#2cite\endcsname##1{\expandafter
     \ifx\csname#2tag##1\endcsname\relax?\else{\rm\csname#2tag##1\endcsname}\fi
    \expandafter\ifx\csname#2tr@cetrue\endcsname\relax\else
     \write\cit@tionsout{#2tag ##1 cited on page \folio.}\fi}%
   \expandafter\gdef\csname#2page\endcsname##1{\expandafter
     \ifx\csname#2page##1\endcsname\relax?\else\csname#2page##1\endcsname\fi
     \expandafter\ifx\csname#2tr@cetrue\endcsname\relax\else
     \write\cit@tionsout{#2tag ##1 cited on page \folio.}\fi}%
   \expandafter\gdef\csname#2tag\endcsname##1{\global\advance
     \csname#2tagno\endcsname by 1%
   \expandafter\ifx\csname#2check##1\endcsname\relax\else%
\fi
   \expandafter\xdef\csname#2check##1\endcsname{}%
   \expandafter\xdef\csname#2tag##1\endcsname
     {#1\number\csname#2tagno\endcsnameX}%
   \write\t@gsout{#2tag ##1 assigned number \csname#2tag##1\endcsname\space
      on page \number\count0.}%
   \csname#2tag##1\endcsname}}%
\def\m@kecs #1tag #2 assigned number #3 on page #4.%
   {\expandafter\gdef\csname#1tag#2\endcsname{#3}
   \expandafter\gdef\csname#1page#2\endcsname{#4}}
\def\re@der{\ifeof\t@gsin\let\next=\relax\else
    \read\t@gsin to\t@gline\ifx\t@gline\v@idline\else
    \expandafter\m@kecs \t@gline\fi\let \next=\re@der\fi\next}
\def\t@gs#1{\def\alltgs{}\m@ketag#1e\m@ketag#1s\m@ketag\t@d@l p
    \m@ketag\gl@b@l r \n@wread\t@gsin\openin\t@gsin=\jobname.tgs \re@der
    \closein\t@gsin\n@wwrite\t@gsout\openout\t@gsout=\jobname.tgs }
\outer\def\localtags{\t@gs\l@c@l}
\outer\def\globaltags{\t@gs\gl@b@l}
\outer\def\newlocaltag#1{\m@ketag\l@c@l{#1}}
\outer\def\newglobaltag#1{\m@ketag\gl@b@l{#1}}

\def\t@gsoff#1,{\def\@{#1}\ifx\@\empty\let\next=\relax\else\let\next=\t@gsoff
   \expandafter\gdef\csname#1cite\endcsname{\relax}
   \expandafter\gdef\csname#1page\endcsname##1{?}
   \expandafter\gdef\csname#1tag\endcsname{\relax}\fi\next}
\def\verbatimtags{\let\ift@gs=\iffalse\ifx\alltgs\relax\else
   \expandafter\t@gsoff\alltgs,\fi}
\catcode`\X=11 \catcode`\@=\active
\localtags
%
\setbox200\hbox{$\scriptscriptstyle \data $}
\global\newcount\numpunt
\hoffset=0.cm
\baselineskip=14pt  
\parindent=12pt
\lineskip=4pt\lineskiplimit=0.1pt
\parskip=0.1pt plus1pt

\hyphenation{small}
\newcount\mgnf\newcount\tipi\newcount\tipoformule\newcount\greco

\tipi=2          
\tipoformule=0   


\global\newcount\numsec
\global\newcount\numfor
\global\newcount\numtheo
\global\advance\numtheo by 1

\def\senondefinito#1{\expandafter\ifx\csname#1\endcsname\relax}

\def\SIA #1,#2,#3 {\senondefinito{#1#2}%
\expandafter\xdef\csname #1#2\endcsname{#3}\else
\write16{???? ma #1,#2 e' gia' stato definito !!!!} \fi}

\def\etichetta(#1){(\veroparagrafo.\veraformula)%
\SIA e,#1,(\veroparagrafo.\veraformula) %
\global\advance\numfor by 1%
\write15{\string\FU (#1){\equ(#1)}}%
\write16{ EQ #1 ==> \equ(#1) }}

\def\letichetta(#1){\veroparagrafo.\verotheo
\SIA e,#1,{\veroparagrafo.\verotheo}
\global\advance\numtheo by 1
\write15{\string\FU (#1){\equ(#1)}}
\write16{ Sta \equ(#1) == #1 }}

\def\tetichetta(#1){\veroparagrafo.\veraformula 
\SIA e,#1,{(\veroparagrafo.\veraformula)}
\global\advance\numfor by 1
\write15{\string\FU (#1){\equ(#1)}}
\write16{ tag #1 ==> \equ(#1)}}

\def\FU(#1)#2{\SIA fu,#1,#2 }

\def\etichettaa(#1){(A\veroparagrafo.\veraformula)%
\SIA e,#1,(A\veroparagrafo.\veraformula) %
\global\advance\numfor by 1%
\write15{\string\FU (#1){\equ(#1)}}%
\write16{ EQ #1 ==> \equ(#1) }}

\def\BOZZA{
\def\alato(##1){%
 {\rlap{\kern-\hsize\kern-1.4truecm{$\scriptstyle##1$}}}}%
\def\aolado(##1){%
 {
{
 \rlap{\kern-1.4truecm{$\scriptstyle##1$}}}}} 
}

\def\alato(#1){}
\def\aolado(#1){}

\def\veroparagrafo{\number\numsec}
\def\veraformula{\number\numfor}
\def\verotheo{\number\numtheo}

\def\Eq(#1){\eqno{\etichetta(#1)\alato(#1)}}
\def\eq(#1){\etichetta(#1)\alato(#1)}
\def\leq(#1){\leqno{\aolado(#1)\etichetta(#1)}}
\def\teq(#1){\tag{\aolado(#1)\tetichetta(#1)\alato(#1)}}
\def\Eqa(#1){\eqno{\etichettaa(#1)\alato(#1)}}
\def\eqa(#1){\etichettaa(#1)\alato(#1)}
\def\eqv(#1){\senondefinito{fu#1}$\clubsuit$#1
\write16{#1 non e' (ancora) definito}%
\else\csname fu#1\endcsname\fi}
\def\equ(#1){\senondefinito{e#1}\eqv(#1)\else\csname e#1\endcsname\fi}

\def\Lemma(#1){\aolado(#1)Lemma \letichetta(#1)}%
\def\Theorem(#1){{\aolado(#1)Theorem \letichetta(#1)}}%
\def\Proposition(#1){\aolado(#1){Proposition \letichetta(#1)}}%
\def\Corollary(#1){{\aolado(#1)Corollary \letichetta(#1)}}%
\def\Remark(#1){{\noindent\aolado(#1){\bf Remark \letichetta(#1).}}}%
\def\Definition(#1){{\noindent\aolado(#1){\bf Definition 
\letichetta(#1)$\!\!$\hskip-1.6truemm}}}
\def\Example(#1){\aolado(#1) Example \letichetta(#1)$\!\!$\hskip-1.6truemm}

\def\include#1{
\openin13=#1.aux \ifeof13 \relax \else
\input #1.aux \closein13 \fi}

\openin14=\jobname.aux \ifeof14 \relax \else
\input \jobname.aux \closein14 \fi
\openout15=\jobname.aux


{\count255=\time\divide\count255 by 60 \xdef\hourmin{\number\count255}
        \multiply\count255 by-60\advance\count255 by\time
   \xdef\hourmin{\hourmin:\ifnum\count255<10 0\fi\the\count255}}

\def\oramin{\hourmin }

\def\data{\number\day/\ifcase\month\or january \or february \or march \or april
\or may \or june \or july \or august \or september
\or october \or november \or december \fi/\number\year;\ \oramin}

\newcount\pgn \pgn=1
\def\foglio{\number\numsec:\number\pgn
\global\advance\pgn by 1}
\def\foglioa{A\number\numsec:\number\pgn
\global\advance\pgn by 1}

\footline={\rlap{\hbox{\copy200}}\hss\tenrm\folio\hss}

\def\TIPIO{
\font\setterm=amr7 
\def \settepunti{\def\rm{\fam0\setterm}
\textfont0=\setterm   
\normalbaselineskip=9pt\normalbaselines\rm }\let\nota=\settepunti}

\def\TIPITOT{
\font\twelverm=cmr12
\font\twelvei=cmmi12
\font\twelvesy=cmsy10 scaled\magstep1
\font\twelveex=cmex10 scaled\magstep1
\font\twelveit=cmti12
\font\twelvett=cmtt12
\font\twelvebf=cmbx12
\font\twelvesl=cmsl12
\font\ninerm=cmr9
\font\ninesy=cmsy9
\font\eightrm=cmr8
\font\eighti=cmmi8
\font\eightsy=cmsy8
\font\eightbf=cmbx8
\font\eighttt=cmtt8
\font\eightsl=cmsl8
\font\eightit=cmti8
\font\sixrm=cmr6
\font\sixbf=cmbx6
\font\sixi=cmmi6
\font\sixsy=cmsy6
\font\twelvetruecmr=cmr10 scaled\magstep1
\font\twelvetruecmsy=cmsy10 scaled\magstep1
\font\tentruecmr=cmr10
\font\tentruecmsy=cmsy10
\font\eighttruecmr=cmr8
\font\eighttruecmsy=cmsy8
\font\seventruecmr=cmr7
\font\seventruecmsy=cmsy7
\font\sixtruecmr=cmr6
\font\sixtruecmsy=cmsy6
\font\fivetruecmr=cmr5
\font\fivetruecmsy=cmsy5
\textfont\truecmr=\tentruecmr
\scriptfont\truecmr=\seventruecmr
\scriptscriptfont\truecmr=\fivetruecmr
\textfont\truecmsy=\tentruecmsy
\scriptfont\truecmsy=\seventruecmsy
\scriptscriptfont\truecmr=\fivetruecmr
\scriptscriptfont\truecmsy=\fivetruecmsy
\def \eightpoint{\def\rm{\fam0\eightrm}
\textfont0=\eightrm \scriptfont0=\sixrm \scriptscriptfont0=\fiverm
\textfont1=\eighti \scriptfont1=\sixi   \scriptscriptfont1=\fivei
\textfont2=\eightsy \scriptfont2=\sixsy   \scriptscriptfont2=\fivesy
\textfont3=\tenex \scriptfont3=\tenex   \scriptscriptfont3=\tenex
\textfont\itfam=\eightit  \def\it{\fam\itfam\eightit}%
\textfont\slfam=\eightsl  \def\sl{\fam\slfam\eightsl}%
\textfont\ttfam=\eighttt  \def\tt{\fam\ttfam\eighttt}%
\textfont\bffam=\eightbf  \scriptfont\bffam=\sixbf
\scriptscriptfont\bffam=\fivebf  \def\bf{\fam\bffam\eightbf}%
\tt \ttglue=.5em plus.25em minus.15em
\setbox\strutbox=\hbox{\vrule height7pt depth2pt width0pt}%
\normalbaselineskip=9pt
\let\sc=\sixrm  \let\big=\eightbig  \normalbaselines\rm
\textfont\truecmr=\eighttruecmr
\scriptfont\truecmr=\sixtruecmr
\scriptscriptfont\truecmr=\fivetruecmr
\textfont\truecmsy=\eighttruecmsy
\scriptfont\truecmsy=\sixtruecmsy }\let\nota=\eightpoint}

\newfam\msbfam   
\newfam\truecmr  
\newfam\truecmsy 
\newskip\ttglue
\ifnum\tipi=0\TIPIO \else\ifnum\tipi=1 \TIPI\else \TIPITOT\fi\fi

\def\e{\epsilon}

\def\E{{I\kern-.25em{E}}}
\def\N{{I\kern-.25em{N}}}
\def\M{{I\kern-.25em{M}}}
\def\R{{I\kern-.25em{R}}}
\def\Z{{Z\kern-.425em{Z}}}
\def\1{{1\kern-.25em\hbox{\rm I}}}
\def\eu{{1\kern-.25em\hbox{\sm I}}}

\def\C{{I\kern-.64em{C}}}
\def\P{{I\kern-.25em{P}}}
\def\eop{{ \vrule height7pt width7pt depth0pt}\par\bigskip}

\def \ba {\beta_a}


\def\chap #1#2{\line{\ch #1\hfill}\numsec=#2\numfor=1}

\def\ba{{\backslash}}

\def\sqr#1#2{{\vcenter{\vbox{\hrule height.#2pt
     \hbox{\vrule width.#2pt height#1pt \kern#1pt
   \vrule width.#2pt}\hrule height.#2pt}}}}
\def\qed{ $\mathchoice\sqr64\sqr64\sqr{2.1}3\sqr{1.5}3$} 


\newcount\foot
\foot=1
\def\note#1{\footnote{${}^{\number\foot}$}{\ftn #1}\advance\foot by 1}
\def\tag #1{\eqno{\hbox{\rm(#1)}}}
\def\frac#1#2{{#1\over #2}}

\def\text#1{\quad{\hbox{#1}}\quad}
\def\newpage{\vfill\eject}

\def\thanks{\noindent{\bf Aknowledgements: }}



\font\ch=cmbx12

\font\ftn=cmr8

\font\it=cmti10
\font\sm=cmr7

%
\catcode`\X=12\catcode`\@=11
\def\n@wcount{\alloc@0\count\countdef\insc@unt}
\def\n@wwrite{\alloc@7\write\chardef\sixt@@n}
\def\n@wread{\alloc@6\read\chardef\sixt@@n}
\def\crossrefs#1{\ifx\alltgs#1\let\tr@ce=\alltgs\else\def\tr@ce{#1,}\fi
   \n@wwrite\cit@tionsout\openout\cit@tionsout=\jobname.cit 
   \write\cit@tionsout{\tr@ce}\expandafter\setfl@gs\tr@ce,}
\def\setfl@gs#1,{\def\@{#1}\ifx\@\empty\let\next=\relax
   \else\let\next=\setfl@gs\expandafter\xdef
   \csname#1tr@cetrue\endcsname{}\fi\next}
\newcount\sectno\sectno=0\newcount\subsectno\subsectno=0\def\r@s@t{\relax}
\def\resetall{\global\advance\sectno by 1\subsectno=0
  \gdef\firstpart{\number\sectno}\r@s@t}
\def\resetsub{\global\advance\subsectno by 1
   \gdef\firstpart{\number\sectno.\number\subsectno}\r@s@t}
\def\v@idline{\par}\def\firstpart{\number\sectno}
\def\l@c@l#1X{\firstpart.#1}\def\gl@b@l#1X{#1}\def\t@d@l#1X{{}}
\def\m@ketag#1#2{\expandafter\n@wcount\csname#2tagno\endcsname
     \csname#2tagno\endcsname=0\let\tail=\alltgs\xdef\alltgs{\tail#2,}%
  \ifx#1\l@c@l\let\tail=\r@s@t\xdef\r@s@t{\csname#2tagno\endcsname=0\tail}\fi
   \expandafter\gdef\csname#2cite\endcsname##1{\expandafter
     \ifx\csname#2tag##1\endcsname\relax?\else{\rm\csname#2tag##1\endcsname}\fi
    \expandafter\ifx\csname#2tr@cetrue\endcsname\relax\else
     \write\cit@tionsout{#2tag ##1 cited on page \folio.}\fi}%
   \expandafter\gdef\csname#2page\endcsname##1{\expandafter
     \ifx\csname#2page##1\endcsname\relax?\else\csname#2page##1\endcsname\fi
     \expandafter\ifx\csname#2tr@cetrue\endcsname\relax\else
     \write\cit@tionsout{#2tag ##1 cited on page \folio.}\fi}%
   \expandafter\gdef\csname#2tag\endcsname##1{\global\advance
     \csname#2tagno\endcsname by 1%
   \expandafter\ifx\csname#2check##1\endcsname\relax\else%
\fi
   \expandafter\xdef\csname#2check##1\endcsname{}%
   \expandafter\xdef\csname#2tag##1\endcsname
     {#1\number\csname#2tagno\endcsnameX}%
   \write\t@gsout{#2tag ##1 assigned number \csname#2tag##1\endcsname\space
      on page \number\count0.}%
   \csname#2tag##1\endcsname}}%
\def\m@kecs #1tag #2 assigned number #3 on page #4.%
   {\expandafter\gdef\csname#1tag#2\endcsname{#3}
   \expandafter\gdef\csname#1page#2\endcsname{#4}}
\def\re@der{\ifeof\t@gsin\let\next=\relax\else
    \read\t@gsin to\t@gline\ifx\t@gline\v@idline\else
    \expandafter\m@kecs \t@gline\fi\let \next=\re@der\fi\next}
\def\t@gs#1{\def\alltgs{}\m@ketag#1e\m@ketag#1s\m@ketag\t@d@l p
    \m@ketag\gl@b@l r \n@wread\t@gsin\openin\t@gsin=\jobname.tgs \re@der
    \closein\t@gsin\n@wwrite\t@gsout\openout\t@gsout=\jobname.tgs }
\outer\def\localtags{\t@gs\l@c@l}
\outer\def\globaltags{\t@gs\gl@b@l}
\outer\def\newlocaltag#1{\m@ketag\l@c@l{#1}}
\outer\def\newglobaltag#1{\m@ketag\gl@b@l{#1}}

\def\t@gsoff#1,{\def\@{#1}\ifx\@\empty\let\next=\relax\else\let\next=\t@gsoff
   \expandafter\gdef\csname#1cite\endcsname{\relax}
   \expandafter\gdef\csname#1page\endcsname##1{?}
   \expandafter\gdef\csname#1tag\endcsname{\relax}\fi\next}
\def\verbatimtags{\let\ift@gs=\iffalse\ifx\alltgs\relax\else
   \expandafter\t@gsoff\alltgs,\fi}
\catcode`\X=11 \catcode`\@=\active
\localtags
%
\setbox200\hbox{$\scriptscriptstyle \data $}
\global\newcount\numpunt
\magnification=\magstep1
\hoffset=0.cm
\baselineskip=12pt  
\parindent=12pt
\lineskip=4pt\lineskiplimit=0.1pt
\parskip=0.1pt plus1pt

\hyphenation{small}

\count0=1

\def\qed{\lower 2 pt \vbox{\hrule width 8 pt height 8 pt depth 0 pt}}
\def\boxin#1{\lower 3.5 pt \vbox{\hrule \hbox{\strut \vrule{} #1 \vrule} 
\hrule}}

 2

 \count0=1
 
\centerline{\bf A Sharp analog of Young's Inequality on $S^N$}
\centerline{\bf and Related Entropy Inequalities}

\def\be{{\vec e}}
\def\bu{{\vec u}}
\def\bv{{\vec v}}
\def\bw{{\vec w}}
\def\bz{{\vec z}}
\def\ba{{\vec a}}
\def\bb{{\vec b}}
\def\bt{{\vec t}}
\def\bc{{\vec c}}
\def\by{{\vec y}}
\overfullrule -0pt

\vskip 1.5truecm
 \centerline{E. A. Carlen\footnote{$^1$}{\eightpoint Work partially supported by U.S. National Science Foundation 
 grant DMS 03-00349.}
\quad
 E. H. Lieb\footnote{$^2$}{\eightpoint Work partially supported by U.S. National Science Foundation grant PHY-0139984.}\quad
M. Loss$^1$}
\vfootnote{}{\baselineskip = 12pt{\eightpoint\copyright 2004 by the authors.
Reproduction of this
article, in its entirety, by any means is permitted for non-commercial
purposes.}}

\medskip
\centerline{1. School of Mathematics, Georgia Tech, Atlanta GA 30332}
\medskip
\centerline{2. Departments of Mathematics and Physics, Jadwin Hall,}
\centerline{Princeton University, P.O. Box 708, Princeton NJ 08544}

\bigskip
\bigskip {\baselineskip = 12pt\narrower{\noindent {\bf Abstract }\/
We prove a sharp analog of Young's 
inequality on $S^N$, and deduce from it certain sharp entropy inequalities.
The proof turns on  constructing a nonlinear heat flow that drives trial functions
to optimizers in a monotonic manner.
This strategy also works for the generalization of Young's inequality on $R^N$ 
to more than three functions, and leads to significant new information about 
the optimizers and the constants.}}

\medskip
\noindent{ Math reviews Classification Numbers:  43A15, 52A40, 82C40}


\newpage
 \chap {1. Introduction }1 

\bigskip

\def\br{{\vec p}}

This paper concerns further generalizations of the generalized Young's inequality due
to Brascamp and Lieb [\rcite{BL}], which we now recall.

For any $N\ge M$, let $\ba_1,\ba_2,\dots,\ba_N$ be any set of $N$ non zero vectors spanning
$\R^M$. Let $f_1,f_2,\dots,f_N$ be any set on $N$ non negative measurable functions on $\R$.
Given numbers $p_j$ with $1 \le p_j \le \infty$
for $j=1,2\dots,N$,
form the vector
$$\br = (1/p_1,1/p_2,\dots,1/p_N)\ ,\Eq(rdef)$$
and define 
$$D(\br) = \sup\left\{{
\int_{\R^M}\prod_{j=1}^N f_j(\ba_j\cdot x){\rm d}^Nx \over
\prod_{j=1}^N\|f_j\|_{p_j}}\ :\ f_j \in L^{p_j}(\R) \qquad j =1,2,\dots, N\ \right\}\ .\Eq(bl1)$$ 

A theorem in [\rcite{BL}] reduces the computation of $D(\br)$ to a 
{\it finite dimensional} variational problem:
Let ${\cal G}$ denote the set of all centered Gaussian
functions on $\R$. That is, $g(x) \in {\cal G}$ if and only if $g(x) = ce^{-(sx)^2/2}$ for some $s >0$
and some constant $c$. Define $D_{{\cal G}}(\br)$ by
$$D_{{\cal G}}(\br)  = \sup\left\{{
\int_{\R^M}\prod_{j=1}^N g_j(\ba_j\cdot x){\rm d}^Nx \over
\prod_{j=1}^N\|g_j\|_{p_j}}\ :\ g_j \in {\cal G}\qquad j =1,2,\dots, N\ \right\}\ .\Eq(bl)$$ 
It is proved in [\rcite{BL}] that $D(\br) = D_{{\cal G}}(\br)$, and hence
$$\int_{\R^M}\prod_{j=1}^N f_j(\ba_j\cdot x){\rm d}^Nx \le  D_{{\cal G}}(\br)
\prod_{j=1}^N\|f_j\|_{p_j}\ ,\Eq(ineqform)$$
This can be used to explicit  compute sharp constants in certain cases. For instance,
when $M=2$ and $N=3$, $D_{{\cal G}}(\br)$ may be evaluated, and this gives
the sharp constant in the classical Young's inequality for convolutions.\footnote{*}{\eightpoint
The sharp constant in Young's inequality for convolutions was obtained by Beckner at the same time that Brascamp 
and Lieb obtained their 
more general result. Beckner's results  do not address the case of more than three functions, 
which is the focus here.}

The first part of this paper  concerns  a version of 
this generalized Young's inequality for functions on the sphere $S^{N-1}$. 
Our generalization was motivated by statistical 
mechanical considerations, and was devised to prove a sharp entropy inequality for 
probability denisities on $S^{N-1}$ 
which is also presented below.
There are by now several alternative proofs of the Brascamp--Lieb inequality 
for functions on $\R^N$ (e.g., [\rcite{L}], [\rcite{Ba}] and [\rcite{Ba2}]), but none of these
seem to be readily adaptable to the consideration of functions on $S^{N-1}$, and it was necessary to develop a new approach. 

The new approach, it turns out, leads to a very simple 
proof of the original theorem in 
$\R^N$, and enables us to strengthen the  original theorem in several respects, 
clearing up some questions left open by the authors cited above. In particular, we resolve a conjecture of Barthe 
whose incisive work in  [\rcite{Ba2}] settled many of the questions about non negative optimizers for the 
variational problem \eqv(bl1). We also obtain new information on the constants. For any given choice of 
$\{\ba_1,\dots,\ba_N\}$, we give an explicit formula for the (unique) choice of $\br$ that minimizes $D(\br)$, 
as well as the minimum value, which we refer to as the ``best best constant'' in the generalized Young's inequality.

We shall proceed to these results along the path which led to them, and begin
by recalling some facts that motivated the investigation of a spherical analog of \eqv(ineqform).

Let $\mu_N$ denote the  uniform Borel probability measure on
$S^{N-1}(\sqrt{N})$, the  sphere of radius $\sqrt{N}$ in $\R^N$, and let
Let $\gamma_N$ denote the Gaussian probability measure
$${\rm d}\gamma_N = (2\pi)^{-N/2}e^{-|\vec v|^2/2}{\rm d}^Nv\ . \Eq(gau)$$
We can consider $\mu_N$ as a probability measure on $\R^N$, supported on $S^{N-1}(\sqrt{N})$, and then it 
is a familiar fact that for large values of $N$, ${\rm d}\gamma_N \approx {\rm d}\mu_N$.
In the considerations that motivated our investigation, a vector 
$$\bv =(v_1,v_2,\dots,v_N)\Eq(vecnot)$$ in $\R^N$ represents the velocities of $N$ one dimensional particles.
Under any sort of evolution of the particle system that conserves kinetic energy, ${\displaystyle \sum_{j=1}^Nv_j^2}$
will be constant. Supposing that its initial value is $N$, at any later time 
the state of the system will be given by a point in $S^{N-1}(\sqrt{N})$. The uniform probability measure
$\mu_N$ is called the {\it microcanonical ensemble} in statistical mechanics. The proability measure
${\rm d}\gamma_N$ on the other hand would be called the {\it canonical ensemble} for this system. 
The principle of {\it equivalence of ensembles} is a cornerstone of equilibrium statistical mechanics.
For this simple system, it reduces to the statement that for any fixed positive integer $k$, and
any  bounded measurable function $\phi(v_1,v_2,\dots,v_k)$
of the first $k$ velocities only,
$$\lim_{N\to \infty}\left(\int_{S^{N-1}(\sqrt{N})}\phi(v_1,v_2,\dots,v_k){\rm d}\mu_N -
\int_{\R^N}\phi(v_1,v_2,\dots,v_k){\rm d}\gamma_N \right) = 0\ .$$

However, the equivalence of ensembles goes only so far. 
A fundamental qualitative difference between $\gamma_N$ and $\mu_N$ is that under the
first probability measure, the coordinate functions are independent random variables, 
while under the second they are not.  This lack of independence has an important quantitative effect 
that {\it does not diminish with increasing} $N$ if we consider functions of all of the velocities at once,
as we now explain.

Before going further, it will be convenient to make a change of scale, and consider the unit sphere.
The factors of $\sqrt{N}$ that are necessary for comparison to the Gaussian measure ${\rm d}\mu_N$
will not be helpful in the next paragraphs. Therefore, let $\mu$ denote the  uniform Borel probability
measure on $S^{N-1}$, the unit sphere in $\R^N$.
Let $\pi_j$ be the $j$th coordinate function. That is, 
$$\pi_j(\bv) = v_j \in [-1,1]. \Eq(kpidef)$$
Consider functions $f_j$ defined on the interval $[-1,1]$ and pull them back to the sphere
via the coordinate function $\pi_j$, i.e., $f_j(\pi_j(\vec v))$. 
We denote this function also by $f_j$. It will be clear from the context which of these
functions is meant.

Because
$\sum_{j=1}^N \pi_j(\vec v)^2 = 1$,
the coordinate functions are not independent random variables, and hence,
given $N$ functions $f_j$ on $[-1,1]$, the quantities
$$\int_{S^{N-1}}\left(\prod_{j=1}^N f_j \right){\rm d}\mu  
\qquad{\rm and}\qquad \prod_{j=1}^N\left(\int_{S^{N-1}} f_j{\rm d}\mu\right)\ .\Eq(print3)$$
need not be equal.
In fact, simple examples show that it is possible for the integral on the left in \eqv(print3) 
to diverge while all of the integrals on the right
are finite. However, according to the following theorem, such a divergence
is not possible if each $f_j$ is square integrable. Indeed, the product of the $L^2$ 
norms of the $f_j$ controls the integral of $\prod_{j=1}^Nf_j$ 
in the strongest way that one could hope.
In what follows, 
$\|\cdot\|_{L^p(S^{N-1})}$ will denote an $L^p$ norm with respect to $\mu$ on $S^{N-1}$.

\bigskip
\noindent{\bf Theorem 1} {\it For all $N\ge 2$, given non--negative  measurable functions
$f_j$,
$j= 1,2,\dots,N$, on $[-1,1]$, 
$$\int_{S^{N-1}}\left(\prod_{j=1}^N f_j(v_j)\right){\rm d}\mu \le
\prod_{j=1}^N \|f_j\|_{L^p(S^{N-1})}\ .\Eq(bound)$$
for all $p \geq 2$. Moreover, for each $p<2$,
there exist functions $f_j$  so that 
$\|f_j\|_{L^p(S^{N-1})} < \infty$ for each $j$, while the integral on the left side of
\eqv(bound) diverges. Finally, for every $p \geq 2$ and $N\ge3$, there is equality in \eqv(bound) if and only if each $f_j$ is constant.}
\bigskip

It is natural to refer to \eqv(bound) as a spherical version of the generalized Young's inequality
\eqv(ineqform). The resemblance is accentuated if we write $v_j = \be_j\cdot\bv$ where the $\be_j$ are the
standard basis vectors in $\R^N$. The proof that we  give for Theorem 1 can be adapted to prove
a generalization in which other vectors other than the $\be_j$ are considered, but this is 
not needed for the application that we now describe.

The inequality \eqv(bound) implies a sharp entropy inequality for probability densities $F$ on $S^{N-1}$. 
Indeed, let $F$ be any probability density on $S^{N-1}$, and then, for each $j=1,2,\dots,N$, let $f_j$
denote  the conditional expectation 
of $F$ given $v_j$. This is a linear operation, and we define the operator $P_j$ by
$$ P_jF = f_j\ .$$

In more analytic terms, $f_j$ is the function on $[-1,1]$ so that for 
all bounded measurable functions $\phi$
on  $[-1,1]$,
$$\int_{S^{N-1}}\phi(v_j)F(\bv){\rm d}\mu = \int_{S^{N-1}}\phi(v_j)f_j(v_j){\rm d}\mu\ .$$
As is evident from the definition,
for square integrable $F$,  $f_j = P_jF$ is just the orthogonal projection in $L^2(S^{N-1})$  of $F$
onto the subspace consisting of functions depending only on $v_j$; i.e., 
measurable with respect to the sigma algebra generated by $\pi_j$.

There is yet another relation worth bearing in mind. To explain, 
introduce the one dimensional marginal $\nu_N$ of $\mu$:
For $N\ge 3$, and any function $\phi$ on $[-1,1]$,
$$\int_{S^{N-1}}\phi(v_1){\rm d}\mu = \int_{[-1,1]}\phi(v){\rm d}\nu_N$$
where
$${\rm d}\nu_N = {|S^{N-2}|\over |S^{N-3}|}(1 - v^2)^{(N-3)/2}{\rm d}v\ .\Eq(marg)$$
Here, $|S^{m-1}|$ denotes the surface area of the $m-1$ dimensional unit sphere in $\R^m$; and $|S^0| = 2$.  Then, 
$f_j(v){\rm d}\nu_N$ is the marginal distribution of $v_j$ under $F(\bv){\rm d}\mu$.
Whenever we refer to the $j$th marginal $f_j$ of a probability density $F$ on $S^{N-1}$, we shall mean that
$f_j$ is related to $F$ in exactly this manner.

Clearly, each of the $f_j$ is a probability density on $S^{N-1}$. For any probability density $F$ on $S^{N-1}$
the entropy of $F$ is $S(F)$ defined by
$$S(F) = \int_{S^{N-1}}F \ln F{\rm d}\mu\ ,\Eq(entdef)$$
and likewise, the entropy of the marginal is given by
$$S(f_j) = \int_{[-1,1]}f_j \ln f_j{\rm d}\nu_N =
 \int_{S^{N-1}}f_j\ln f_j{\rm d}\mu\ .$$
How do the entropies of the marginals $f_j$   compare with the entropy of their parent density $F$?
The following theorem provides an answer:

\medskip 
\noindent{\bf Theorem 2} {\it For all $N\ge 2$, given any probability density $F$ on $S^{N-1}$, let $f_j$
be the $j$th marginal of $F$ for $j=1,2,\dots,N$. Then
$$\sum_{j=1}^N S(f_j) \le 2S(F)\ ,\Eq(subadd)$$
and the constant $2$ on the right side of \eqv(subadd) is the best possible.}
\medskip

The inequality \eqv(subadd) may be compared to the familiar 
{\it subadditivity of the entropy} inequality on $\R^N$:
Let $G$ be any probability density on $\R^N$ with respect to 
${\rm d}\gamma_N$, and let $g_j$ denote its $j$th marginal, which is obtained by integrating out
all of the variables except $v_j$. In this case, there is no relation among the coordinate functions.  Hence
$$\int_{\R^N}\prod_{j=1}^N g_j {\rm d}\gamma_N = \prod_{j=1}^N\left(\int_{\R^N} g_j{\rm d}\gamma_N \right) = 1\ ,$$
so that $H = \prod_{j=1}^N g_j$ is another probability density on $\R^N$. Then by Jensen's inequality,
$$\eqalign{
0 &\le \int_{\R^N}\left({G\over H}\right)\ln \left({G\over H}\right) H{\rm d}\gamma_N\cr
&= \int_{\R^N} G \ln G {\rm d}\gamma_N - \int_{\R^N} G \ln  H{\rm d}\gamma_N\cr
&= \int_{\R^N} G \ln G {\rm d}\gamma_N - \sum_{j=1}^N\int_{\R^N} G \ln g_j{\rm d}\gamma_N\ ,\cr}$$
and there is equality if and only if $G = H$.
Defining the entropy $S(G)$ of a density $G$ realtive to ${\rm d}\gamma_N$ by 
${\displaystyle S(G) = \int_{\R^N} G \ln G {\rm d}\gamma_N}$,
this says
$$\sum_{j=1}^NS(g_j) \le S(G)\Eq(subaddg)$$
with equality if and only if $G = H$. Note the difference between \eqv(subaddg) and \eqv(subadd):
The latter requires an extra factor of $2$ on the right, independent of $N$. This is due to the
dependence of the coordinate functions $v_j$ resulting from the constraint 
$\sum_{j=1}^N v_j^2 = 1$.

The difference between  \eqv(subaddg) and \eqv(subadd) is especially striking given 
the close relation between ${\rm d}\mu$ and ${\rm d}\gamma_N$.
The inequality in Theorem 2 does not depend on the radius of the sphere,
since the uniform measure is normalized, and so Theorem 2 says that there is a {\it dimension 
independent} departure from the equivalence of ensembles as measured by subadditivity of the entropy.

This dimension independence would not be guessed by linearizing the inequality in Theorem 2 about
$F=1$; {\it it is a non--perturbative effect}. The natural perturbative calculation would suggest
that the difference between   \eqv(subaddg) and \eqv(subadd) ``washes out'' with increasing $N$, 
as we now explain.

Consider a probability denisty $F$ on $S^{N-1}$ of the form
$$F=1+\varepsilon H$$
where $H$ is bounded and orthogonal to $1$ in  $L^2(S^{N-1})$. Then
$f _j = P_jF = 1 + \varepsilon h_j$
where $h_j = P_jH$. Of course $h_j$ is also orthogonal to  $1$ in  $L^2(S^{N-1})$.

A simple and frequently encountered computation gives us
$$S(F) = {\varepsilon^2\over 2}\|H\|^2_{L^2(S^{N-1})} + {\cal O}(\varepsilon^3)\qquad
{\rm and}\qquad
\sum_{j=1}^N S(f_j) = {\varepsilon^2\over 2}\|h_j\|^2_{L^2(S^{N-1})}+ {\cal O}(\varepsilon^3)\ .$$
Since 
$$\|h_j\|^2_{L^2(S^{N-1})} = \langle P_j H, P_j H \rangle_{L^2(S^{N-1})}
 = \langle  H, P_j H \rangle_{L^2(S^{N-1})}\ ,$$
if we define the operator 
$$P = {1\over N}\sum_{j=1}^NP_j\ ,$$
we have that
$${\sum_{j=1}^NS(P_jF) \over S(F)} = N{ \langle  H, P H \rangle_{L^2(S^{N-1})}
\over \|H\|^2_{L^2(S^{N-1})}  } + {\cal O}(\varepsilon)\ .$$

An optimist might then hope that the supremum of $\sum_{j=1}^N S(P_jF)/S(F)$ taken over all probability
densities $F$ would be given by $C_N$ where
$$C_N = \sup\left\{N{ \langle  H, P H \rangle_{L^2(S^{N-1})}
\over \|H\|^2_{L^2(S^{N-1})}  }\ :\ H\in L^2(S^{N-1})\ , \ \langle H,1\rangle_{L^2(S^{N-1})}= 0\right\}\ .\Eq(kac)$$
The computation of the supremum is an eigenvalue problem, and has been done in 
[\rcite{CCL1}; Theorems 1.2 and 2.1]. The result is
$$
C_N= 1+{3 \over N+1} \ .
$$
The surplus over $1$, namely $3/(N+1)$, measures the  ``departure from independence'' 
as a function of $N$.
Thus, if one considers densities $F$ that deviate only slightly from the uniform density,
one gets a correction term to the constant $1$ in the  Gaussian entropy inequality \eqv(subaddg) 
that ``remembers'' the
dependence of the coordinates on the sphere, but which vanishes as $N \to \infty$. 

The precise size
of this ``departure from independence'' as a function of $N$ is crucial in some problems of 
non--equilbrium statisitical mechanics.  The computation of $C_N$ was at the
heart of recent progess in computing the rate of relaxation to equilibrium in  kinetic theory 
by direct consideration of an $N$ body system, as proposed long ago by Mark Kac. 
For more details, see [\rcite{CCL1}] and [\rcite{CCL2}].  

The fact that for more general densities $F$, the  correction term to the Gaussian entropy inequality \eqv(subaddg)
does not vanish as $N \to \infty$ complicates the estimation of  rates of realxation in entropic terms 
for $N$ body systems in kinetic theory.  This said, we turn to the proof of Theorem 2.

\medskip
\noindent{\bf Proof of Theorem 2:} Let $F$ be any probability density on $S^{N-1}$, $N\ge 2$, and let $f_j$, $j=1,2,\dots,N$
be its marginals. Then since $f_j$ is a probability density, $\|f_j^{1/2}\|_{L^2(S^{N-1})} =1$. As a consequence of
Theorem 1, if we define $C$ by
$$C = \int_{S^{N-1}}\left(\prod_{j=1}^N f^{1/2}_j\right){\rm d}\mu\ ,$$
we have $C \le 1$. 

Suppose that $C =0$. Then, $\prod_{j=1}^N f_j = 0$ almost everywhere, and so $\sum_{j=1}^N\ln f_j = -\infty$
almost everywhere. This would imply
$$- \infty = \int_{S^{N-1}}F\left(\sum_{j=1}^N\ln f_j\right){\rm d}\mu = \sum_{j=1}^NS(f_j)\ .$$
This is impossible, since by Jensen's inequality, $S(f_j) \ge 0$ for each $j$. Therefore, $0 < C< 1$,
and we may define a probability density $H$ on $S^{N-1}$ through
$$H = {1\over C}\prod_{j=1}^N f^{1/2}_j\ .$$
As above, we now apply Jensen's inequality to conclude that
$$0 \le \int_{S^{N-1}}\left({F\over H}\right)\ln \left({F\over H}\right) H{\rm d}\mu\ .$$
The right and side is easily seen to be
$$\eqalign{
&\int_{S^{N-1}}F \ln F{\rm d}\mu - \int_{S^{N-1}}F \ln H{\rm d}\mu =\cr
&\int_{S^{N-1}}F \ln F{\rm d}\mu - \sum_{j=1}^N \int_{S^{N-1}}F \ln f_j^{1/2}{\rm d}\mu + \ln(C) =\cr
&\int_{S^{N-1}}F \ln F{\rm d}\mu - {1\over 2}\sum_{j=1}^N \int_{S^{N-1}}f_j \ln f_j{\rm d}\mu + \ln(C) \ .\cr}$$
Since $\ln(C) <1$ unless each $f_j =1$, the inequality is proved, with equality holding only when
each $f_j =1$. 
The fact that the constant cannot be less than $2$ in the inequality follows by finding a trial function 
that we present in the Appendix. \eop
\medskip

As discussed above, the factor of $2$ in Theorem 2 is a 
correction to the classical subadditivity of the entropy that is required
on account of the dependence of the coordinate functions due to the constraint $\sum_{j=1}^Nv_j^2 =1$. 
The  remarkable fact that the size of this effect is independent
of $N$ depends on the specific nature of the constraint, and is not a general fact.

For example, consider the planar constraint
$\sum_{j=1}^N v_j = 0$
on $\R^N$, and  let $P_{N-1}$ denote  the hyperplane specified by this constraint. Let $\tilde \mu$ be 
a centered,
isotropic Gaussian probability measure on $P_{N-1}$. As we explain below, a special case of the 
Brascamp--Lieb Theorem yields the sharp inequalitiy
$$\int_{P_{N-1}}\left(\prod_{j=1}^N f(v_j)\right)\tilde \mu \le \prod_{j=1}^N\| f_j\|_{L^{N/(N-1)}(P_{N-1})}\ .\Eq(boundan)$$

This is an analog of \eqv(bound) for the planar constraint.  Notice however, that this time the $L^p$
indices depend on $N$, and diminish towards $1$ as $N$ increases.

Just as Theorem 2 follows from Theorem 1,  one obtains an entropy subadditivity inequality for the planar 
constraint from \eqv(boundan).
Given a probability density $F$ on $P_{N-1}$ with respect to the reference measure $\tilde \mu$, define 
the marginal densities as above. Then the analog of \eqv(subadd) is the inequality
$$\sum_{j=1}^N S(f_j) \le {N\over N-1}S(F)\ .\Eq(subaddan)$$
This time,  since $\lim_{N\to \infty}N/(N-1) =1$, the effect of the contraint,
as far as subaddititvity of the entropy is concerned,
diminishes to zero as $N$ tends to infinity.

The connection between \eqv(boundan) and  Young's inequality is revealing.  To see the connection, we  change of variables.
Let  $\be_j$, $j =1,2,\dots, N$ be the standard basis vectors in $\R^N$. 
Let $\bu_j$ be the normalized orthogonal projection of $\be_j$
onto the hyperplane $P_{N-1}$. One easily checks that  for $i \ne j$,
$$\bu_i \cdot \bu_j = -{1\over N-1}$$
and that for $\bv$ in $P_{N-1}$,
$$v_i = \vec v \cdot \be_i = \sqrt{N-1\over N}\bv\cdot \bu_j\Eq(rel1)$$
and 
$$\sum_{j=1}^N|\bu \cdot \bv|^2 = {N\over N-1}|\bv|^2\ .\Eq(rel2)$$
For convenient constants, choose a scale so that the Gaussian denisty is
$M(\bv) = e^{-\pi |\bv|^2}$
Defining the single variable funtions
${\displaystyle g_j(y) = f_j\left(\sqrt{N-1\over N}y\right)e^{- \pi (N-1)y^2/N}}$,
we have from \eqv(rel1) and \eqv(rel2) that
$$\prod_{j=1}^N f_j(v_j)M(\bv) = \prod_{j=1}^N g_j(\bu_j \cdot \bv)\ ,\Eq(rel3)$$
and if ${\rm d}{\cal L}$ denotes Lebesgue measure on $P_{N-1}$,
$$\int_{P_{N-1}}\left(\prod_{j=1}^N f_j(v_j)\right){\rm d}\tilde \mu = 
\int_{P_{N-1}}\left(\prod_{j=1}^N g_j(\bu_j\cdot\bv)\right){\rm d}{\cal L}\ .\Eq(rel4)$$
Furthermore, for each $j$,
$\int_{P_{N-1}}|f_j(v_j)|^{N/(N-1)}{\rm d}\tilde \mu = \int_{\R}|g_j(y)|^{N/(N-1)}{\rm d}y$.
Therefore, identifying $P_{N-1}$ with $\R^{N-1}$ in the obvious way, \eqv(boundan) is equivalent to the inequality
$$\int_{\R^{N-1}}\left(\prod_{j=1}^N g_j(\bu_j\cdot\bv)\right){\rm d}{\cal L} \le
\prod_{j=1}^N\|g_j\|_{L^{N/(N-1)}(\R)}\ ,\Eq(rel6)$$
which is a special case of the Brascamp--Lieb generalization of Young's inequality. 

In fact, for $N=3$, it follows from 
the sharp form of the classical Young's inequality for convolutions. To see this, let $s = \bu_1\cdot \bv$ and 
$t = \bu_3\cdot \bv$,
and notice that since $\bu_1 + \bu_2 + \bu_3 =0$, we have that $-(s+t) = \bu_2\cdot\bv$. A simple computation reveals that
${\rm d}{\cal L} = (2/\sqrt{3}){\rm d}s{\rm d}t$,
so that the $N=3$ case of \eqv(rel6) becomes
$$\int_{\R^2}g_1(s)g_2(-s-t)g_3(t){\rm d}s{\rm d}t \le 
{\sqrt{3}\over 2}\|g_1\|_{L^{3/2}(\R)}\|g_2\|_{L^{3/2}(\R)}\|g_3\|_{L^{3/2}(\R)}\ .$$
This in turn is equivalent to the inequality
$$\|g_1 * g_2\|_{L^{3}(\R)} \le {\sqrt{3}\over 2}\|g_1\|_{L^{3/2}(\R)}\|g_2\|_{L^{3/2}(\R)}\ ,$$
which is sharp.  We now turn to the proof of Theorem 1.

\bigskip
\chap {2. Proof of Young's inequality on $S^{N-1}$ }2 
\medskip

We prove Theorem 1 using a non-linear heat semigroup. For $1\le i,j \le N$, $i\ne j$
let 
$$L_{i,j}  = v_i{\partial \over \partial v_j} - 
v_j{\partial \over \partial v_i}\ .$$
The Laplacian on $S^{N-1}$ is the operator
$$\Delta  = \sum_{i<j}L_{i,j}^2 = {1\over 2}\sum_{i\ne j}L_{i,j}^2\ .\Eq(lap)$$
The normalization of the gradient on $S^{N-1}$ implicit in this is convenient; 
for smooth functions $f$ and $g$, we write
$$\nabla f\cdot \nabla g = \sum_{i<j}L_{i,j}f L_{i,j}g = {1\over 2}\sum_{i\ne j}L_{i,j}f L_{i,j}g\ ,\Eq(bracket)$$
and $|\nabla f|^2  = \nabla f\cdot \nabla f$.

Now fix any $p\ge 1$. For any smooth, non negative 
function $g$ in $L^p(S^{N-1})$, and any $t>0$,
define
$$g(\bv,t) = \left(e^{t\Delta }g^p(\bv)\right)^{1/p}\ .\Eq(pev)$$
The first thing to observe is that $g(\cdot,t)$ will be smooth and strictly positive for all $t>0$,
and the $L^p(S^{N-1})$ norm of $g$ is conserved under this evolution:
$$\|g(\cdot, t)\|_{L^p(S^{N-1})} = \|g\|_{L^p(S^{N-1})}\Eq(cons)$$
for all $t\ge 0$.

The second thing to observe is that if $g$ depends only on $v_j$ for some $j$,
so does $g(\cdot, t)$. The reason is that  $g$ depends only on $v_j$ if and only if
$g$ is invariant under all rotations that fix the $j$th coordinate axis, and these rotations 
commute with the Laplacian. We write $g(v_j,t)$ to denote the evolution of such a function.

The third thing to observe is that the evolution, though non-linear, has the semigroup property:
For all $s,t>0$,
$$g(\bv,s+t) = \left(e^{s\Delta }g^p(\bv,t)\right)^{1/p}\ .\Eq(semi)$$

The fourth thing to observe is that 
$$\lim_{t\to \infty}g(\bv,t) =  \|g\|_{L^p(S^{N-1})}\Eq(lim)$$
uniformly in $\bv$. 

Finally,
a simple computation shows that for any smooth, positive  function $g$ on $S^{N-1}$,
$${ \partial \over \partial t}g(v,t)\bigg|_{t= 0} = 
{1\over p}g^{1-p}\Delta g^p = \Delta g +  (p-1){|\nabla g|^2\over g}\ .\Eq(pgen)$$

\medskip
\noindent{\bf Lemma 2.1} {\it Consider any $N$ 
non negative  functions
$g_1,g_2,\dots, g_N$ in $L^2([-1,1],{\rm d}\nu_N)$. Use $p=2$ and $g_j$ in place of $g$ in \eqv(pev)
to define $g_j(v_j,t)$. Then by the smooting properties of the heat equation, the function $\phi(t)$
defined by
$$\phi(t) = \int_{S^{N-1}}\prod_{j=1}^N g_j(v_j,t){\rm d}\mu$$
is differentiable for all $t>0$, and is right continuous at $t= 0$. Moreover, introducing
the functions $h_k$ and $G$  defined by
$$h_j(v_j,t) = \ln g_j(v_j,t)\qquad k =1,2\dots, N\qquad{\rm and}\qquad G = \prod_{j=1}^Ng_j\ ,\Eq(hGdef)$$ 
$${{\rm d}\over {\rm d}t}\phi(t) = 
{1\over 2}\int_{S^{N-1}}\sum_{i\ne k}
 \left[ (L_{i,k} h_k) - (L_{i,k} h_i)\right]^2 G{\rm d}\mu\ .\Eq(comp6)$$
}

\medskip

\noindent{\bf Proof} The statements about smoothness and continuity require no justification. 
Taking $p=2$ and $g= g_k(v_k,t)$ in \eqv(pgen),
we have
$${\partial \over \partial t}g_k(v_k,t)
 = \Delta g_k(v_k,t) +  {|\nabla g_k(v_k,t)|^2\over g_k(v_k,t)}\ .$$
Hence, supressing the arguments on the right,
$${{\rm d}\over {\rm d}t}\left(
\int_{S^{N-1}}\prod_{j=1}^N g_j(v_j,t){\rm d}\mu\right)\bigg|_{t=0}
 = \sum_{k=1}^N \int_{S^{N-1}}
 \left(\Delta g_k +  {|\nabla g_k|^2\over g_k}\right)\prod_{\ell =1,\ell\ne k}^Ng_\ell {\rm d}\mu\ .$$
 
The integral on the right can be written as
$$ \int_{S^{N-1}}\sum_{k=1}^N 
\left(\Delta g_k \right)\prod_{\ell =1,\ell\ne k}^Ng_\ell {\rm d}\mu +
\int_{S^{N-1}}\sum_{k=1}^N  
\left({|\nabla g_k|^2\over g_k}\right)\prod_{\ell =1,\ell\ne k}^Ng_\ell {\rm d}\mu\ .\Eq(two)$$
clearly, the second integral on the right is non negative. 
We therefore examine the first integral.

Observe that $L_{i,j}g_k = 0$ unless either $i= k$ or $j=k$. Therefore,
$$\int_{S^{N-1}}\sum_{k=1}^N 
\left(\Delta g_k \right)\prod_{\ell =1,\ell\ne k}^Ng_\ell {\rm d}\mu  = 
\int_{S^{N-1}}\sum_{k=1}^N 
\left(\sum_{i<k}L^2_{i,k} g_k  + 
\sum_{j> k}L^2_{k,j} g_k\right)\prod_{\ell =1,\ell\ne k}^Ng_\ell {\rm d}\mu \Eq(comp1)$$

Integrating by parts,
$$\int_{S^{N-1}}\sum_{k=1}^N \left(\sum_{i<k}L^2_{i,k}g_k\right)
\prod_{\ell =1,\ell\ne k}^Ng_\ell {\rm d}\mu = -
\int_{S^{N-1}}\sum_{k=1}^N \left(\sum_{i<k}L_{i,k}g_k L_{i,k}g_i\right) 
\prod_{\ell =1,\ell \ne i,k}^Ng_\ell {\rm d}\mu \ .\Eq(comp2)$$

Using the notations in \eqv(hGdef), 
the integral on the right  side of
\eqv(comp2) is
$$-\int_{S^{N-1}}\sum_{k=1}^N \sum_{i<k}\left(L_{i,k}h_kL_{i,k}h_i\right)G{\rm d}\mu\ .$$
Doing the same integration by parts on the remaining terms in \eqv(comp1), and substituting
$i$ for $j$, we have
$$
\int_{S^{N-1}}\sum_{k=1}^N 
\left(\Delta g_k \right)\prod_{\ell =1,\ell\ne k}^Ng_\ell {\rm d}\mu  =
- \int_{S^{N-1}}\sum_{i\ne k}\left(L_{i,k}h_kL_{i,k}h_i\right)G{\rm d}\mu\ .\Eq(comp4)$$

With the same notations, we have
$$\eqalign{
\int_{S^{N-1}}\sum_{k=1}^N  
\left({|\nabla g_k|^2\over g_k}\right)\prod_{\ell =1,\ell\ne k}^Ng_\ell {\rm d}\mu
 &= \int_{S^{N-1}}\sum_{k=1}^N(\nabla h_k)^2G{\rm d}\mu\cr
 &= \int_{S^{N-1}}\sum_{i\ne k}(L_{i,k} h_k)^2G{\rm d}\mu\cr
 &= {1\over 2}\int_{S^{N-1}}\sum_{i\ne k}
 \left[ (L_{i,k} h_k)^2 + (L_{i,k} h_i)^2\right] G{\rm d}\mu\ .\cr}\Eq(comp5)$$
Combining \eqv(comp4) and \eqv(comp5) we see that
$${{\rm d}\over {\rm d}t}\left(
\int_{S^{N-1}}\prod_{j=1}^N g_j(v_j,t){\rm d}\mu\right)\bigg|_{t=0} =
{1\over 2}\int_{S^{N-1}}\sum_{i\ne k}
 \left[ (L_{i,k} h_k) - (L_{i,k} h_i)\right]^2 G{\rm d}\mu\ .$$
This is \eqv(comp6). \eop
\medskip

\noindent{\bf Proof of Theorem 1:}\ \  By Lemma 2.1, the difference between the right and left hand sides
of \eqv(bound) is
$$\eqalign{
&\int_0^\infty  \left( {{\rm d}\over {\rm d}t}
\int_{S^{N-1}}\prod_{j=1}^N g_j(v_j,t){\rm d}\mu \right) {\rm d}t =\cr 
{1\over 2}&\int_0^\infty\left(\int_{S^{N-1}}\sum_{i\ne k}
 \left[ (L_{i,k} h_k(v_k,t)) - (L_{i,k} h_i(v_i,t))\right]^2 G{\rm d}\mu\right){\rm d}t \ge 0\ .\cr}$$
This proves the inequality. 

Also, it is now clear that for all $t>0$, each $h_k$ is smooth and bounded, and $G$ is strictly positive, so that
there is equality 
in \eqv(bound) if and only if
$\left[ (L_{i,k} h_k) - (L_{i,k} h_i)\right]^2 = 0$
for all $t>0$, all $\bv$ and all $i\ne k$.

 Fixing $t$, $i$ and $k$, this requires
$$v_i h'_k(v_k) = -v_kh'_i(v_i)\Eq(eeqq)$$
This implies that for some constant $c$, 
$\left[\matrix{h'_i(v_i)\cr h'_k(v_k)\cr}\right] =c\left[\matrix{-v_i\cr \phantom{-}v_k\cr}\right]$ 
for all values of $v_i$ and $v_k$. Hence, for all $i\ne k$,
$h_i'$ and $h_k'$ are linear functions with slopes of the same magnitude but opposite signs.
For $N\ge 3$, the signs of all pairs cannot be opposite 
unless all of the slopes are zero.
This concludes the proof that there is equality in \eqv(bound) if and only if each
of the functions there is constant.\footnote{*}{\eightpoint The analysis of \eqv(eeqq) in a
preprint of this paper contained an error. This was pointed out 
and corrected in a private communication from Shannon Starr, to whom we are grateful.}

In the appendix, there is an explicit example showing 
that \eqv(bound) cannot hold if $L^2(S^{N-1})$ is 
replaced by $L^p(S^{N-1})$ for any $p<2$.  In fact, it is shown that for any $p<2$, there is a function $f$ so that with 
$f_j = f$ for all $j$, the left hand side of \eqv(bound) is infinite, and the
right hand side is finite. 
Alternatively, one can see that if \eqv(bound)
did hold with $2$ repalced by some $p < 2$, then Theorem 2 would hold with $2$ replaced
by this value of $p$, which we have seen is not possible.
\eop

The simple heat flow argument that was used to prove Theorem 1 can be adapted to other situations as well. Indeed, one 
could easily consider inequalities for integrals over $S^{N-1}$ of more general products
${\displaystyle \prod_{j=1}^P f_j(\ba_j\cdot \bv)}$. 
The case considered here was $P=N$ and $\ba_j = \be_j$ because that was what was relevant
for Theorem 2. Further generalizations are possible, and may be interesting.

In the next sections, we exhibit the versatility of the method by  showing that a heat flow interpolation between trial 
functions and Gaussians  can be used to prove the
original Brascamp Lieb inequality on $\R^M$.  Barthe [\rcite{Ba}],[\rcite{Ba2}] has
given a proof of this inequality, together with a dual inverse inequality using an interpolation based on optimal mass 
transport. It was somewhat surprising that one could prove the Brascamp Lieb inequality with such a simple heat flow 
interpolation, and after hearing a report on our work, Barthe
and Cordero--Erausquin [\rcite{BaCE}] added to the surprise by showing that a heat flow interpolation could
be used to prove the inverse dual inequality as well.

\bigskip
\chap {3. The generalized Young's inequality on $\R^{M}$ }3
\medskip

We have introduced this inequality in the introduction, 
and shall use the same notation here.  Recall that
for any $N\ge M$,  $\ba_1,\ba_2,\dots,\ba_N$ is a set of $N$ non zero vectors in
$\R^M$. Let $f_1,f_2,\dots,f_N$ be any set of $N$ non negative measurable functions on $\R$,
and consider the integral
$$\int_{\R^M}\prod_{j=1}^N f_j(\ba_j\cdot x){\rm d}^M x\ .\Eq(yo1N)$$
There are certain natural restrictions on the 
underlying set of vectors $\ba_j$.
First of all,  $\{\ba_1,\ba_2,\dots,\ba_N\}$ must span $\R^M$ for
the integral in \eqv(bl1) to  possibly converge.  
Second, it is natural to assume that no pair of vectors $\ba_i$ and $\ba_j$ are proportional; if
they were, we could combine two factors into one in the integrand \eqv(bl1).   {\it These assumptions will
be in force throughout the following sections}.

As  before, given numbers $p_j$ with $1 \le p_j \le \infty$
for $j=1,2\dots,N$,
form the vector
$\br = (1/p_1,1/p_2,\dots,1/p_N)$, and define 
$D(\br)$ and $D_{\cal G}(\br)$ through  \eqv(bl1) and \eqv(bl) respectively.

The Brascamp and Lieb Theorem asserts that $D(\br) = D_{{\cal G}}(\br)$. As in the proof of Theorem 2, we shall use a 
non linear semigroup based on an appropriately chosen heat kernel
to interpolate between arbitrary trial functions and Gaussian optimizers. The appropriate choice of
the heat kernel depends on both $\br$ and the vectors $\{\ba_1,\dots,\ba_N\}$. We shall show
in this section 
that the desired heat kernel exist whenever the supremum is attained the Gaussian variational problem \eqv(bl) for given  
$\br$ and  $\{\ba_1,\dots,\ba_N\}$.   Note that the supremum being attained means that  there are numbers 
$0 <s_j < \infty$ so that with $g_j(y) = e^{-(s_j^2y^2)/2}$, 
$$\int_{\R^M}\prod_{j=1}^N g_j(\ba_j\cdot x){\rm d}^Nx = D_{\cal G}(\br)
\prod_{j=1}^N\|g_j\|_{p_j}$$
In this case, we shall say that the Gaussian variational problem has an optimizer, and we identify the optimizer 
with the vector in $\R^N$ whose $j$th entry is $s_j$. 

\medskip

\noindent{\bf Theorem 3.1} {\it Let $\{\ba_1,\ba_2,\dots,\ba_N\}$ be a set of vectors
spanning $\R^M$ and suppose that the vector $\br$ is such that the 
Gaussian variational problem \eqv(bl) has a maximizer. Then $D(\br)  = D_{{\cal G}}(\br)$. }
\medskip

By itself, this theorem is contained in  the Brascamp--Lieb Theorem, which asserts that 
 $D(\br)  = D_{{\cal G}}(\br)$ in general. However, as we shall see in the next section, Theorem 3.1
 provides the essential reduction of \eqv(bl1) to \eqv(bl), and to complete the analysis and prove the full result, 
 one needs only certain facts about the Gaussian variational problem. For the most part, the facts we need
 are contained in the work of Barthe [\rcite{Ba2}], so that once we have proved Theorem 3.1, our work is largely done.  
 The rest of this section is devoted to the proof of Theorem 3.1.

\medskip

As preparation for the proof, let $R$ be any invertible $M\times M$ matrix, and consider the heat semigroup $e^{tL}$
generated by
$$L = \nabla \cdot R^tR\nabla\ .\Eq(yo2)$$ For each $j$, and each $t>0$, define $f_j(\cdot,t)$
by
$$f_j(t,\ba_j\cdot x) = \left( e^{tL}f_j^{p_j}(\ba_j\cdot x)\right)^{1/p_j}\ .\Eq(yo4)$$

Since $L$ commutes with translations, the set of  functions on $\R^M$
of the form $f(\ba\cdot x)$  is invariant under $e^{tL}$ . In fact, 
for any bounded function $f_0$ on $\R$, for all $t>0$,
$e^{tL}f_0(\ba\cdot x) = f(t,\ba\cdot x)$
where $f(t,y)$ is the solution of
$${\partial \over \partial t}f(t,y) = 
|R\ba|^2{\partial^2 \over \partial y^2}f(t,y)\qquad f(0,y) = f_0(y)\ .\Eq(yo3)$$
The fundamental solution of \eqv(yo3) is 
${\displaystyle g(t,y) = {1\over \sqrt{4\pi |R\ba|^2 t}}e^{-y^2/(4|R\ba|^2 t)}}$.
Therefore, with pointwise convergence, 
 $$\lim_{t\to\infty}t^{1/2p_j}f_j(t,t^{1/2}y) = 
 {\|f_j\|_{p_j}\over (4\pi |R\ba |^2)^{1/2p_j}}e^{-y^2/(4|R\ba|^2 p_j)}\ .\Eq(ber1)$$
 Let $g_j(y)$ denote the centered Gaussian function defined  by the right hand side
 of \eqv(ber1).  
 Note also that for each $j$ and $t$, 
 $$\|f_j(t,\cdot)\|_{p_j} = \|f_j\|_{p_j}  =   \|g_j\|_{p_j}\ .\Eq(ber2)$$

If we assume that each $f_j$ is bounded and has compact support, then it
 is possible to obtain simple Gaussian bounds on each $f_j(t,y)$ from which, using  \eqv(ber1) and
the obvious  dominated convergence argument,  it follows that 
 $$\lim_{t\to \infty}\int_{\R^M}\prod_{j=1}^Nt^{1/2p_j}f_j(t,t^{1/2}(\ba_j\cdot x)){\rm d}^Mx
  = \int_{\R^M}\prod_{j=1}^N g_j  (\ba_j\cdot x)){\rm d}^Mx\ .$$
 Moreover, by the scale invariance that obtains under \eqv(yo22),
 $$\int_{\R^M}\prod_{j=1}^Nt^{1/2p_j}f_j(t,t^{1/2}(\ba_j\cdot x)){\rm d}^Mx =
 \int_{\R^M}\prod_{j=1}^N f_j(t,(\ba_j\cdot x)){\rm d}^Mx$$
so that
$$\lim_{t\to \infty}\int_{\R^M}\prod_{j=1}^N f_j(t,(\ba_j\cdot x)){\rm d}^Mx
  =  \int_{\R^M}\prod_{j=1}^N g_j  (\ba_j\cdot x)){\rm d}^Mx\ .$$

It now follows that {\it if} we can choose $R$ so that ${\displaystyle
\int_{\R^M}\prod_{j=1}^N f_j(t,(\ba_j\cdot x)){\rm d}^Mx}$ is a non decreasing function of $t$,
then
$${\int_{\R^M}\prod_{j=1}^N f_j(\ba\cdot x){\rm d}^Nx \over
\prod_{j=1}^N\|f_j\|_{p_j}} \le
{\int_{\R^M}\prod_{j=1}^N g_j(\ba\cdot x){\rm d}^Nx \over
\prod_{j=1}^N\|g_j\|_{p_j}}\ .$$

By this argument, proof of the Brascamp--Lieb Theorem is reduced to finding a fixed matrix $R$
so that the function
$\eta(t)$ defined by
$$ \eta(t) = \int_{\R^M}\prod_{j=1}^N f_j(t,\ba_j\cdot x){\rm d}^Nx\ \Eq(ber4)$$
is non--decreasing 
where  $f_j(t,y)$ is determined through the choice of $R$  by \eqv(yo2) and \eqv(yo4).

If this is to work at all, the Gaussian functions $g_j$ defined by the limit in \eqv(ber1) must be maximizers for
the variational problem \eqv(bl1), and certainly for the variational problem \eqv(bl). We can gain insight 
into how $R$ must be chosen by considering the Euler--Lagrange equation for \eqv(bl).

For each $j$,
let $\tilde g_j$ be the centered Gaussian function given by $\tilde g_j(y) = e^{-(s_jy_j)^2/2}$.
Then a simple calculation reveals that
$$2\ln\left({\int_{\R^M}\prod_{j=1}^N \tilde g_j(\ba_j\cdot x){\rm d}^Mx \over
\prod_{j=1}^N\|\tilde g_j\|_{p_j}}\right)
= \sum_{j=1}{1\over p_j}2\ln(s_j) - \ln\left(\det\left(AS^2A^t\right)\right)\ ,\Eq(yo66)$$
where $S$ is the diagonal $N\times N$ matrix whose $j$th diagonal entry is $s_j$, and $A$ is the
$M \times N$ matrix whose $j$th column is $\ba_j$; i.e., $A = [\ba_1,\ba_2,\dots,\ba_N]$. 
(This notation $A = [\ba_1,\ba_2,\dots,\ba_N]$ will be used repeatedly in what follows.)

Introduce the variables $t_1,t_2,\dots,t_N$ by $t_j = \ln(s_j^2)$. Let $T$ be the diagonal matrix whose $j$
diagonal entry is $t_j$. 

Define the function $\phi$ on $\R^N$ by
$$\phi(t_1,t_2,\dots,t_N) = Tr\left(\ln\left(Ae^TA^t\right)\right)\ .\Eq(yo63)$$
Since $\ln\left(\det\left(Ae^TA^t\right)\right) = Tr\left(\ln\left(Ae^TA^t\right)\right)$,
we have from \eqv(bl) and \eqv(yo66) that
$$2\ln(D_{{\cal G}}(\br)) = \sup_{t_1,t_2,\dots,t_N}\left\{\sum_{j=1}{1\over p_j}t_j - 
\phi(t_1,t_2,\dots,t_N)\right\}\Eq(yo65)$$ 

A simple calculation shows that
$$\eqalign{
{\partial\over \partial t_j}\phi(t_1,t_2,\dots,t_N)  &= 
e^{t_j}\ba_j\cdot (Ae^TA^t)^{-1}  \ba_j\cr
&= (s_j\ba_j)\cdot (AS^2A^t)^{-1}(s_j\ba_j)\ .\cr}\Eq(yo64)$$
Therefore, the Euler--Lagrange equation for the 
optimization problem in \eqv(yo65) is 
$$\eqalign{
{1\over p_j} &= s_j\ba_j \cdot(AS(AS)^t)^{-1}s_j\ba_j\cr
&= \be_j \cdot (AS)^t (AS(AS)^t)^{-1} (AS)\be_j\ ,\cr}\Eq(bl2)$$
where $\be_j$ is the $j$th standard basis vector in $\R^N$.  Notice  that since 
${\rm rank}(A) = M$,  the matrix
$(AS)^t (AS(AS)^t)^{-1} (AS)$ is just the orthogonal projection onto the range of $AS$.

We now show that  if the supremum in \eqv(yo65) is a attained, so that there is a positive
diagonal matrix $S$ satisfying \eqv(bl2), then we can choose 
$$R = (AS(AS)^t)^{-1/2}\ ,\Eq(choice)\ ,$$
and with this choice,  the function $\eta(t)$ defined by \eqv(ber4) is non--decreasing. The key is the following lemma:

\medskip

\noindent{\bf Lemma 3.2} {\it Let $f_1,f_2,\dots,f_N$ be $N$ bounded, non--negative measurable functions 
 on $\R$ with compact support.
Let $R$ be any invertible $M\times M$ matrix, and consider the heat semigroup $e^{tL}$
generated by $L = \nabla \cdot R^tR\nabla$. For each $j$, and each $t>0$, define $f_j(\cdot,t)$
by \eqv(yo4)
and define 
Define the function $\eta(t)$ by \eqv(ber4). Then with $h_j(y,t) = \ln f_j(t,y)$ and 
$F(x,t) = \prod_{j=1}^Nf_j(t,\ba_j\cdot x)$,
$\eta(t)$ is differentiable for $t>0$, and 
$${{\rm d}\over {\rm d}t}\eta(t) = 
\int_{\R^M}\left(\sum_{i,j =1}^M 
h'_i(\ba_i\cdot x,t)Q_{i,j}h'_j(\ba_j\cdot x,t)\right)F(x,t){\rm d}^Mx\ ,\Eq(yo21a)$$
where $Q$ is the $M\times M$ matrix with
  $$Q_{i,j} = \delta_{i,j}p_j|R\ba_j|^2 - R\ba_i\cdot R\ba_j\ .\Eq(yo21)$$}

\medskip
\noindent{\bf Proof:} 
By \eqv(yo3) we have that 
$${\partial \over \partial t}f_j(t,\ba_j\cdot x) = Lf_j(t,\ba_j\cdot x)
+ (p_j -1){|R\nabla f_j(t,\ba_j\cdot x)|^2\over f_j(t,\ba_j\cdot x)}$$
and hence
$${{\rm d}\over {\rm d}t}\eta(t) = 
\sum_{j=1}^N\int_{\R^M}\left[
  Lf_j(\ba_j\cdot x)
+ (p_j -1){|R\nabla f_j(\ba_j\cdot x,t)|^2\over f_j(\ba_j\cdot x,t)}\right]\prod_{i\ne j}
 f_i(\ba_i\cdot x,t){\rm d}^M\ .$$
 Let $h_j  = \ln f_j$, and let $F(x) = \prod_{j=1}^Nf_f(\ba_j\cdot x)\ .$
 Then, integrating by parts in the term containing $L$, and suppressing arguments,
 $$\eqalign{
 {{\rm d}\over {\rm d}t}\eta(t) &= 
  \int_{\R^M}\left[ 
 \sum_{j=1}^N(p_j-1)|R\ba_j|^2 |h'_j|^2 - \sum_{i\ne j}(R\ba_i\cdot R\ba_j)
 h'_i h'_j\right]F(x){\rm d}^M x\cr
&=\int_{\R^M}\left[ \sum_{j=1}^N p_j|R\ba_j|^2 |h'_j|^2 - \sum_{i, j}(R\ba_i\cdot R\ba_j)
 h'_i h'_j\right]F(x){\rm d}^M x
  \cr}\Eq(ber111)$$
Using the definition \eqv(yo21), we have \eqv(yo21a).  \eop
\medskip

\noindent{\bf Proof of Theorem 3.1:} 
We apply Lemma 3.2, and must choose $R$ so that $Q$ is a positive matrix.   
By assumption, there is a maximizer for the Brascamp--Lieb variational problem \eqv(bl), or equivalently \eqv(yo65), 
and hence there 
is a positive diagonal matrix $S$ such that  the Euler--Lagrange equation \eqv(bl2) is satisfied  for each $j$. 
In this case with  $R = (AS(AS)^t)^{-1/2}$,  \eqv(yo21) becomes
$$Q = S^{-1}(I - P)S^{-1}\Eq(qfor)$$
where $P = (AS)^t(AS(AS)^t)^{-1}(AS)$ is the orthogonal projection onto the range of $AS$. 
This is certainly non--negative, and hence whenever the Gaussian variational problem \eqv(bl)
has an optimizer,
$D(\br) = D_{{\cal G}}(\br)$. \eop

\medskip

We close this section by expressing $D(\br)$ in terms of $A$ and $S$ when the supremum in \eqv(yo65) is attained. 
In this case, the optimizing  Gaussians $g_j$ are given by the limit in  \eqv(ber1). 
We may assume that $\|f_j\|_{p_j} =1$ for each $j$.  With $R$ given by \eqv(choice),
the Euler--Lagrange equation \eqv(bl2) says $|R\ba_j|^2 = 1/(s_j^2p_j)$, and hence
${\displaystyle g_j(\ba_j\cdot x) =  \left({p_js^2_j\over 4\pi} \right)^{1\over2 p_j} e^{-(s_j\ba_j\cdot x)^2/4}}$.  
Thus,
$$\prod_{j=1}^N g_j(\ba_j\cdot x)  = \prod_{j=1}^N
 \left({p_js^2_j\over 4\pi} \right)^{1\over 2 p_j}e^{-|ASx|^2/4}\ ,\Eq(ber55)$$
and therefore,
$$\eqalign{
D(\br) &= 
\prod_{j=1}^N\left(\left( p_j s_j^2\right)^{1/(2p_j)}\right)
\left({1\over 4\pi }\right)^{M/2}\int_{\R^M}e^{-|ASx|^2/4}{\rm d}^Mx\cr
& = 
\prod_{j=1}^N\left(\left( p_j s_j^2\right)^{1/(2p_j)}\right){\rm det}(AS^2A^t)^{-1/2}\ .\cr}\Eq(ber45)$$
For future use, note that  \eqv(ber55) can be written as
$$\prod_{j=1}^N g_j(\ba_j\cdot x)  = D(\br)\left({1\over \int_{\R^M}e^{-|ASx|^2/4}{\rm d}^Mx}\right)
e^{-|ASx|^2/4}\ .\Eq(ber56)$$

Note that if $S$ satisfies the Euler--Lagrange equation \eqv(bl2). so does $\lambda S$ for any
$\lambda > 0$. Replacing $S$ by $\lambda S$ in \eqv(ber56), and taking $\lambda$ to infinity,
we obtain $ D(\br)\delta_0$ in the limit, where $\delta_0$ is the point mass at the origin.  This will
be used later on.

\bigskip
\chap {4. The Gaussian  optimization problem}4
\medskip

The analysis in the previous section leads very naturally to the following questions:
\medskip
\noindent{$\bullet$} {\it For which values of $\br$ is $D_{{\cal G}}(\br)$ finite?}
\medskip
\noindent{$\bullet$} {\it For which values of $\br$ is there an optimizer for the Gaussian 
variational problem \eqv(bl)?}
\medskip

These questions have been answered by Barthe [\rcite{Ba2}]. 
(In the special case in which every set of
$M$ vectors chosen from among $\{\ba_1,\dots,\ba_N\}$ is a basis, 
this had been done in [\rcite{BL}]).
The answers may appear unfavorable for our program, since it turns out that in general there exist $\br$
for which $D_{{\cal G}}(\br)$ is finite, but for which
there is no optimizer for the Gaussian 
variational problem. Hence one additional observation is required to deduce the 
Brascamp Lieb Theorem from Theorem 3.1.  

First, we recall Barthe's answer to the first question, which is pleasingly simple:  
Let $K_A$ denote the
convex hull of the vectors
$\bz$ whose entries are either zero or one, and for which the set
$\{\ba_j \ :\ z_j =1\ \}$ is a basis of $\R^M$.
Barthe has proved [\rcite{Ba2}] that $D_{{\cal G}}(\br)$ finite if and only if $\br \in K_A$. 

Note that $K_A$ lies in the hyperplane in $\R^N$ given by the equation $\sum_{j=1}^Nz_j = M$. 
Let $K_A^\circ$
denote the interior of $K_A$ relative to this hyperplane.   
Barthe has also proved in [\rcite{Ba2}]  that 
when $\br\in K_A^\circ$, the supremum in the Gaussian optimization problem \eqv(bl) 
is attained. 

In this section, we give another proof of these results. We do this for two reasons. 
First, our proof gives an alternate characterization of $K_A$ that is directly checkable. 
Second, our proof is based on a {\it partial scale invariance} property of the functional
that we seek to optimize. This  partial scale invariance property of the functional is expressed in
the identity \eqv(factor2) below.
As we shall see,  it   completely 
determines the nature of $K_A$, and it provides a crucial handle on the variational problem 
in case $\br$ is on the boundary of $K_A$.

The obvious  scale invariance argument shows that for $D(\br)$ or even $D_{\cal G}(\br)$
to be finite, it is necessary that
$$\sum_{j=1}^N{1\over p_j} = M\  ,\Eq(yo22)$$
and of course that 
$1 \le  p_j \le  \infty$ for each $j$. 
Indeed, let $\lambda$ be any positive number, and  replace each $f_j(y)$ 
in   \eqv(bl1)  by $f_j(\lambda y)$.
The numerator in   \eqv(bl1) is proportional to $\lambda^{-\sum_{j
=1}^N1/p_j}$, while the denominator  is proportional to $\lambda ^{-M}$. 
This excludes a finite maximum unless \eqv(yo22) holds.   

A somewhat less obvious {\it partial scaling argument} leads to further restrictions on $\br$. 
This depends on a simple identity that is crucial in what follows:
\medskip
\noindent{\bf Lemma 4.1} {\it  Let $\{\ba_1,\ba_2,\dots,\ba_N\}$ be a set of vectors
spanning $\R^N$. Let $S$ be any proper, non--empty subset of $ \{1,2,\dots,N\}$, and let 
let $r = {\rm dim}({\rm span}(\{\ba_j\ :\ j\in S\}))$   Then there are explicitly computable  sets of vectors  
$\{\bb_j\ :\ j\in S\}$ and  
$\{\bc_j\ :\ j\in S^c\}$ such that for any set of non negative functions $f_j$, each bounded and with compact support,
$$\eqalign{
&\int_{\R^M}\prod_{j=1}^N f_j(\ba_j\cdot x){\rm d}^M x  =\cr
&\int_{\R^r}\prod_{j\in S} f_j(\bb_j\cdot y) \left(
\int_{\R^{M-r}}\prod_{j\in S^c} f_j(\bb_j\cdot y + \bc_j\cdot z){\rm d}^{M-r}z\right){\rm d}^r y\ .\cr}\Eq(factor)$$
}
\medskip

\noindent{\bf Proof:} Let $\{\bu_1,\dots,\bu_r\}$ be an orthonormal basis for  the span of
$\{\ba_j\ :\ j\in S\}$. 
Let $\{\bv_{r+1},\dots,\bv_M\}$ be an orthonormal basis for  the orthogonal complement.  Choose the sign of $\bv_M$ so that 
$\det([\bu_1,\dots \bu_r,\bv_{r+1},\dots,\bv_M]) =1$.  Let $U = [\bu_1,\dots \bu_r]$, and let
$V = [\bv_{r+1},\dots,\bv_N]$. Define 
$\bb_j = U^t\ba_j$ and $ \bc_j = V^t\ba_j$.
Likewise define
$y = U^tx$ and $ z = V^t x$.
Then $\ba_j \cdot x = \bb_j\cdot \by + \bc_j\cdot \bz$,
and for $j\in S$, $\ba_j \cdot x = \bb_j\cdot \by $. 
Since ${\rm d}^M x = {\rm d}^r y{\rm d}^{M-r}z$, 
and since, by  construction, $\bc_j = 0$ for $j \in S$, \eqv(factor) follows immediately. \eop
\medskip

To apply this, we rescale in ${\rm span}(\{\ba_j\ :\ j\in S\})$ alone: For each $j\in S$,
replace $f_j$ by $f^{(\lambda)}_j$ where $f^{(\lambda)}_j(y) = 
\lambda^{1/p_j}f_j(\lambda y)$. 
Then $\|f^{(\lambda)}_j\|_{p_j} = \|f_j\|_{p_j}$, so 
that this replacement does not affect the denominator in \eqv(bl1).  Then:
$$\eqalign{
&\int_{\R^M}\prod_{j\in S} f_j^{(\lambda)}(\ba_j\cdot x)
\prod_{j\in S^c} f_j(\ba_j\cdot x){\rm d}^M x =\cr
&\left[\lambda^{-r}\prod_{j\in S}\lambda^{1/p_j}\right]\int_{\R^r}\prod_{j\in S} 
f_j(\bb_j\cdot y) \left(
\int_{\R^{M-r}}\prod_{j\in S^c} f_j(\lambda^{-1}\bb_j\cdot y + \bc_j\cdot z){\rm d}^{M-r} 
z\right){\rm d}^r y\ .\cr}\Eq(factor2)$$

We see that if 
$$\sum_{j\in S}{1\over p_j} >  r = {\rm dim}({\rm span}(\{\ba_j\ :\ j\in S\}))\ ,$$
then the integral in  \eqv(factor2) diverges as $\lambda$ tends to $+\infty$. Since  
$\|f^{(\lambda)}_j\|_{p_j} = \|f_j\|_{p_j}$, 
this means that $D(\br)$ is infinite in this case. These considerations justify the following definitions.
\footnote{*}{\eightpoint
Note that as $\lambda$ tends to zero,
$\prod_{j\in S^c} f_j(\lambda^{-1}\bb_j\cdot y + \bc_j\cdot z)$ tends to zero, and will even vanish idenitically for $\lambda$ 
large enough when the $f_j$ have compact support.  Hence, \eqv(factor2) does not give us information on the relation between 
$\sum_{j\in S} 1/p_j$ and $r(S)$ in the limit as $\lambda$ tends to zero.}
\medskip

\noindent{\bf Definition} Let  $\{\ba_1,\ba_2,\dots,\ba_N\}$ be a given set of vectors spanning $\R^M$. 
For each subset $S$ of $\{1,2\dots,N\}$, define  
$$r(S) =    {\rm dim}({\rm span}(\{\ba_j\ :\ j\in S\}))\ .\Eq(raa)$$
Let $K_A$ denote the subset of $\R^N$ consisiting of vectors $\bz$ such that
$\sum_{j=1}^N z_j = 1$, $0\le z_j \le 1$ for each $j$, and finally
$$\sum_{j\in S}z_j \le r(S)\ .\Eq(kaa)$$  
Define $K_A^\circ$ to be the subset of $K_A$ consisting of vectors $\bz$ satisfying
$$\sum_{j\in S}z_j < r(S)\Eq(kaaint)$$
for all proper, non--empty subsets $S$ of  $\{1,2\dots,N\}$. For later use, we say that a subset $S$
is {\it critical} at $\bz$ if $\sum_{j \in S} z_j= r(S)$ and {\it subcritical}
at $\bz$ if $\sum_{j \in S} z_j < r(S)$.
\medskip

It may seem that we are being inconsistent in our notation, as we have already used $K_A$ to denote a certain
convex hull in our description of Barthe's results. We shall show below that in fact the two sets coincide. 
For present puroses, this is not important, and the definition of $K_A$ shall be the one made just above.

What we have just seen shows that $\br\in K_A$ is a necessary condidtion for $D(\br) < \infty$,
or even $D_{{\cal G}}(\br) < \infty$.   It turns out that it is also sufficient. 

 \medskip
\noindent{\bf Theorem 4.2 (Barthe)} {\it  Let   $\{\ba_1,\ba_2,\dots,\ba_N\}$ be any 
a spanning set of vectors in $\R^M$. Then $D_{\cal G}(\br) < \infty$ if and only if $\br\in K_A$.
Moreover, 
if $\bz \in K_A^\circ $, then  the supremum is attained in the Gaussian 
variational problem \eqv(bl). }
\medskip

Barthe's proof is based on the convex hull description of $K_A$, as mentioned above. At the end of this section we 
give an alternate proof, and show directly that Barthe's convex hull definition of  $K_A$  yields the same set as 
does our definition.  First, we deduce the Brascamp--Lieb Theorem from Theorems 3.1, 4.2 and Lemma 4.1.

 \medskip
\noindent{\bf Theorem 4.3 (Brascamp--Lieb)} {\it  
Let   $\{\ba_1,\ba_2,\dots,\ba_N\}$ be any 
a spanning set of vectors in $\R^M$. Then for all $\br$,  $D_{\cal G}(\br) = D(\br)$.}
\medskip

\noindent{\bf Proof:}    If $\br\in K_A^\circ$, everything is clear. By Theorem 4.2, the Gaussian problem has optimizers, 
and then by Theorem 3.1,  $D_{\cal G}(\br) = D(\br)$.

Therefore, suppose that
$\br \in K_A$, but not $K_A^\circ$.    Then there exists a
non--empty proper subset $S$ of the indices  that is critical; i.e, such that
$\sum_{j\in S}1/p_j = r(S)$. We further take $S$ to have the {\it least cardinality} among all such 
sets.

To apply the identity \eqv(factor), consider 
$$ D_S = \sup\left\{{
\int_{\R^r}\prod_{j\in S} f_j(\bb_j\cdot y){\rm d}^r y \over
\prod_{j\in S}\|f_j\|_{p_j}}\ :\ f_j \in L^{p_j}(\R) \qquad j =1,2,\dots, N\ \right\}\ \Eq(ber201)$$
and
$$D_{S^c} =  \sup\left\{{
\int_{\R^{M-r}}\prod_{j\in S^c} f_j(\bc_j\cdot z){\rm d}^{M-r} z \over
\prod_{j\in S^c}\|f_j\|_{p_j}}\ :\ f_j \in L^{p_j}(\R) \qquad j =1,2,\dots, N\ \right\}\ .\Eq(ber202)$$
Here, as in \eqv(factor), $r= r(S)$.
Clearly, \eqv(factor) yields the bound $D(\br) \le D_SD_{S^c}$.

To obtain the opposite inequality, note that the
the scaling identity $\sum_{j\in S}(1/p_j) = r$ is satisfied, and by the choice of a critical set of minimal 
cardinality,  there are {\it no}
critical subsets of $S$ for the variational problem of computing $D_S$.

Therefore, there is a solution of the Euler--Lagrange equation \eqv(bl2) for \eqv(ber201),
and hence it has Gaussian maximizers. From \eqv(ber56) we see that we can take
these maximizers $g_j$ so that  $\prod_{j\in S} f_j(\bb_j\cdot y)$ is an arbitrarily good approximation of $D_S$ 
times a Dirac mass at the origin. This will eliminate the terms involving $y$ in the second integral in \eqv(factor). 
Hence for any $\epsilon>0$, one can choose the functions  $f_j$, $j\in S$, to be Gaussian
and have
$$
{\int_{\R^r}\prod_{j=1}^N f_j(\bb_j\cdot y){\rm d}^r y \over
\prod_{j=1}^N\|f_j\|_{p_j}}  \ge D_S  \left(
{\int_{\R^{M-r}}\prod_{j\in S^c} f_j(\bc_j\cdot z){\rm d}^{M-r }z \over
\prod_{j\in S^c}\|f_j\|_{p_j}}\right) - \epsilon\ .$$

We are now reduced to proving that the variational problem for $D_{S^c}$
has Gaussian maximizers.  If there are no critical subsets of $S^c$ for this problem, we are done by
Theorem 4.1. Otherwise, ``peel off'' another critical subset. This procedure clearly reduces the cardinality 
of $S^c$ each time, and hence it terminates with a full set of Gaussian trial functions that
come arbitrarily close to the supremum.   This yields the 
identity  $D(\br) = D_S D_{S^c}$ and completes the proof of Theorem 4.3. \eop

\medskip

The remainder of this section is devoted to the proof of Theorem 4.2. We first show that the two definitions of $K_A$
do indeed define the same set.

 \medskip
\noindent{\bf Theorem 4.4} {\it  For any  spanning set of vectors $\{\ba_1,\ba_2,\dots,\ba_N\}$,  $K_A$ is a the
convex hull of the vectors
$\bz$ whose entries are either zero or one, and for which the set
$\{\ba_j \ :\ z_j =1\ \}$ is a basis of $\R^M$.}
\medskip

First we prove a lemma.
 \medskip
\noindent{\bf Lemma  4.5} {\it  Consider any spanning set of vectors $\{\ba_1,\ba_2,\dots,\ba_N\}$, 
and any $\bz$ in $K_A$.  Let $T$ be any non--empty subset of
the indices $\{1,2,\dots,N\}$. 
If there is any  set $S$ of indices containing $T$ that is critical at $\bz$, then there is a least such set $S_0$: 
That is, there is a set $S_0$ containing $T$ that is critical at $\bz$, such that
if $\tilde S$ is any other set containing $T$ that is critical at $\bz$, then $S_0 \subset \tilde S$.
.}
\medskip

\noindent{\bf Proof:} Without loss of generality, we may suppose that  there
is a set of indices containing $T$ that is critical at $\bz$.  Let $S$ be such a set of least cardinality, and let $\tilde S$ 
be any other set of indices containing $T$ that is critical at $\bz$.
Let $V$ be the span of $\{\ba_j\ :\  j \in S\}$, and let  $W$ be the span of $\{\ba_j\ :\  j \in \tilde S\}$.
Clearly,
$$\{\ba_j\ :\  j \in S\} \cap \{\ba_j\ :\  j \in \tilde S\} \subset V\cap W$$
and 
$$\{\ba_j\ :\  j \in S\}\cup \{\ba_j\ :\  j \in \tilde S\} \subset V\cup W\ .$$
{}From the relation
${\rm dim}(V\cap W) +  {\rm dim}(V\cup W)  = {\rm dim}(V) + {\rm dim}(W)$
it then follows that 
$$r(S\cap \tilde S) + r(S\cup \tilde S) \le r(S) + r(\tilde S)\ .\Eq(ber91)$$
Since both $S$ and $\tilde S$ are critical at $z$,   $r(S) + r(\tilde S) = \sum_{j\in S}z_j +
 \sum_{j\in \tilde S}z_j $.  Since $\bz \in K_A$,  $\sum_{j\in S\cup\tilde S}z_j \le r(S\cup\tilde S) $
 thus from \eqv(ber91)
 $$\sum_{j\in S\cup\tilde S}z_j  +  r(S\cap \tilde S)  \le 
  \sum_{j\in S}z_j +
 \sum_{j\in \tilde S}z_j \ .$$
 This implies that
 ${\displaystyle r(S\cap \tilde S)  \le 
  \sum_{j\in S\cap \tilde S}z_j}$
  and since $T$ is non--empty and is a subset of both $S$ and $\tilde S$, $S\cap \tilde S$ is not empty. Hence
  $S\cap \tilde S$ is critical at $z$.  If  $S\cap \tilde S$ were a proper subset of $S$,
  then we would have found a critical subset of strictly smaller cardinality, contrary to the assumption on $S$. 
  Hence $S = S\cap \tilde S$, and so $S\subset \tilde S$. The set 
  $S$ is the set $S_0$ that we seek. \eop
  
 \medskip
\noindent{\bf Proof of Theorem 4.4:}   Suppose that $\bz \in K_A$, and for some $k$, 
$0 < z_k  < 1$. We shall show that in this case, $\bz$ is not extreme.

First, consider the case in which  no critical set contains any indices $j$ for which  $0 < z_j  < 1$.

Since $\sum_{j=1}^Nz_j = M$,  which  is an integer, it must be the case that for some $\ell\ne k$, 
$0 < z_\ell  < 1$.     Since neither $k$ nor $\ell$ belongs to any critical set,  increase (resp. decrease) 
$z_k$ a little, while decreasing (resp. increasing) 
$z_\ell$ a little in such a way that $z_k + z_\ell$ is constant, and 
the increases do not produce any supercritical sets.  Clearly in this case, 
$\bz$ is not extreme.

Second, if there are critical sets containing  indices $j$ for which  $0 < z_j  < 1$, choose one, $S$, of least cardinality.
Since $S$ is critical, $\sum_{j\in S}z_j$ is an integer, and there must be two indices $k$ and $\ell$ in $S$
such that $0 < z_k,z_\ell < 1$. By the lemma and the definition of $S$, $S$ is the smallest critical set 
containing either $k$ or $\ell$.

Clearly,  we can increase $z_k$ a little bit, and decrease $z_\ell$ a little bit
without affecting $z_k + z_\ell$, and hence  without affecting 
$\sum_{j\in S}z_j$. Moreover, the increase in $z_k$ does not increase the value of
$\sum_{j\in \tilde S}z_j$ for any other critical set $\tilde S$ that contains $k$. This is because
$S \subset \tilde S$ by Lemma 4.5 and the definition of $S$, and hence $\tilde S$ also contains $\ell$. \eop

\medskip
\noindent{\bf Proof of Theorem 4.2:}  
Let $\phi_A(t_1,t_2,\dots,t_N) = \ln\left(\det(Ae^TA^t)\right)$. The function $\phi_A$ 
was shown to be convex on $\R^N$ by Brascamp and Lieb. Let 
 $\phi_A^*$ denote its Legendre transform:
 $$\phi_A^*(z_1,z_2,\dots,z_N) = 
  \sup_{t_1,t_2,\dots,t_N}\left\{\sum_{j=1}z_jt_j - 
\phi_A(t_1,t_2,\dots,t_N)\right\}\Eq(yo65b)$$ 
By \eqv(yo65), determining the set of vectors $\br $ for which $D_{{\cal G}}(\br) < \infty$ is 
the same as determining the set of vectors $\bz$ for which $\phi_A^*(\bz) < \infty$.

Next, recall a formula of Brascamp and Lieb, which can be deduced from the Cauchy--Binet formula:
$$\det\left(Ae^TA^t\right) = \sum_{|S| = M} t_S \det \left( A_S A_S^t\right)\ ,\Eq(yo181)$$
where $t_S = \exp\left(\sum_{j\in S}t_j\right)$.
Here, we use the following notation: If 
$S = \{j_1,j_2,\dots, j_k\}$,
 $A_S$ denotes the $M \times k$ matrix
$$A_S = [\ba_{j_1},\ba_{j_2},\dots,\ba_{j_k}]\Eq(yo85)$$
As shown in [\rcite{BL}], the convexity of $\phi_A(\bt)$ follows by differentiating $\phi(\bt + r\bv)$ twice with respect
to $r$ using the Schwarz inequality. Here $\bv$ is an arbitrary fixed vector. 

Having made these remarks,     
we first show that $\phi^*(\bz) = \infty$ unless 
$\sum_{j=1}^Nz_j = M$, $0\le z_j \le 1$ for each $j$.

For any constant $c$ and any $\bt$ in $\R^N$, let $\bt_c$ denote the vector in $\R^N$
whose $j$th component is $t_j + c$. From the definition \eqv(yo63), it follows that
$\phi_A(\bt_c) = Mc + \phi_A(\bt)$.
Therefore,
$$\bz\cdot \bt_c - \phi_A(\bt_c) = \left(\sum_{j=1}^Nz_j - M\right)c + \bz\cdot \bt - \phi_A(\bt) \Eq(invar)$$
so that the domain of $\phi_A^*$ lies in the hyperplane
${\displaystyle\sum_{j=1}^Nz_j = M}$.
Further, it follows from \eqv(yo64) that 
${\displaystyle 0 \le {\partial \over \partial t_j}\phi_A(\bt)\le 1}$
since this quantity is is the $j$th diagonal entry of an orthogonal projection. Hence,  every
$\bz$ in $K_A$ is such that $0 \le z_j \le 1$ for each $j$.

Recall the terminology that a subset $S$
is {\it critical} at $\bz$ if $\sum_{j \in S} z_j= r(S)$ and {\it subcritical}
at $\bz$ if $\sum_{j \in S} z_j < r(S)$.

We now show that $\phi^*_A(\bz) < \infty $ if $ \bz \in K_A^\circ $. 
First, note that if  $\sum_{j=1}^Nz_j = M$, then \eqv(invar) reduces to 
$$\bz\cdot \bt_c - \phi_A(\bt_c) =   \bz\cdot \bt - \phi_A(\bt) \ .\Eq(invar2)$$
Hence  in \eqv(yo65b), we may restrict our focus to vectors $\bt$
satisfying $\min_{j=1}^N t_j = 0$. 

For any $\bt = (t_1,\dots,t_N)$, let $\bt^* = (t^*_1,\dots,t^*_N)$ be its decreasing rearrangement. 
By the invariance noted above, we may assume that $t^*_N =0$. 
Let $\pi$ be any permutation so that 
$$t^*_j = t_{\pi(j)}\qquad{\rm for\ all}\quad 1 \le j \le N\ .$$

Let $\tilde S$ be the indices of the pivotal columns in $\pi(A) = [\ba_\pi(1),\ba_\pi(2),\dots,\ba_\pi(N)]$. 
That is, the columns of $A_{\tilde S}$ are the columns in $\pi(A)$ that are not in the span of the columns 
to their left in $\pi(A)$. Since the dimension of the space spanned by the vectors $\ba_1 , \dots \ba_N$
is $M$ we have that $|\tilde S|=M$. By monotonicty of the logarithm and \eqv(yo181),
$$\phi_A(t_1,t_2,\dots,t_N) = \ln(\det(Ae^TA^t)) \ge \sum_{j\in \tilde S}t_j + 
\ln(\det(A_{\tilde S}A_{\tilde S}^t)\ ,\Eq(blkey)$$
and hence it suffices to find a lower bound on   
$$\sum_{j\in \tilde S}t_j  - \sum_{j=1}^Nz_j t_j \ . \Eq(tobound) $$

Setting $a_k = 1$ if $\pi^{-1}(k)\in S$, and $a_k = 0$ otherwise and $b_k = z_{\pi(k)}$, \eqv(tobound)
can be written as
$$\sum_{k=1}^N(a_k - b_k)t^*_k\ .$$
The point about this notation is that the vector $(a_1, \dots, a_M, 0, \dots, 0)$ 
which has $N$ elements
strictly majorizes
the vector $(b_1, \dots b_N)$, i.e., 
$$
\sum_{j=1}^k a_k > \sum_{j=1}^k b_k  , k=1, \dots, N-1 \ , \Eq(major1)
$$
and
$$
\sum_{j=1}^N a_k = \sum_{j=1}^N b_k = M \ . \Eq(major2)
$$
The equation \eqv(major2) follows from the definition of the $a_k$'s and $b_k$'s, 
the fact that $|\tilde S|=M$ and
the fact that $\sum_{j=1}^N z_j =M$. The equation \eqv(major1) follows from the definition
of the $a_k$'s and $b_k$'s and the fact that $\bz \in K_A^\circ$, i.e., for every proper
subset $S$ of $\{1, \dots, N\}$, $\sum_{j \in S} z_j < r(S)$. 

Summing by parts, using $t^*_N =0$ and $\sum_{k=1}^Na_k  = \sum_{k=1}^Nb_k$,
$$\sum_{k=1}^N(a_k - b_k)t^*_k = \sum_{k=1}^{N-1}\left(\sum_{j=1}^k(a_j - b_j)\right)(t_k^* - t_{k+1}^*)
\geq c t^*_1 \ ,$$
where $c = \min_{1 \le k \le N-1} \sum_{j=1}^k (a_j-b_j) > 0$.

 Hence
$$\sum_{k\in S}t_k - \sum_{k=1}^Nz_k t_k  \ge  c \max_k(t_k)\ ,$$
which, together with the bound \eqv(blkey), yields
 $$\bz\cdot \bt  - \phi_A(\bt) \le  -\ln(\det(A_{\tilde S}A_{\tilde S}^t) - 
 c \max_k(t_k)\ .$$
 Therefore, as  any of the variables $t_1, \dots, t_{N-1}$ tend to infinity (recall that without loss
 $t_N$ can be chosen to be zero), $\bz\cdot\bt - \phi_A(t_1,\dots,t_N)$
 tends to $-\infty$, and so $\phi^*_A(\bz)  < \infty$. 
 By the convexity of $\phi_A$, proved by Brascamp and Lieb, the
 supremum in \eqv(yo65b) is attained in this case.

 It remains to show that   $\phi^*_A(\bz)  < \infty$ for all $\bz$ in $K_A$.  This is an easy consequence of
 Theorem 4.4. Suppose that  $\br$ is one of the verticies of $K_A$.  If $p_j = \infty$, we may as well replace $f_j$ by $1$
 in \eqv(bl1).  Therefore, there are effectively only $M$ vectors and functions. Letting $S$ denote the 
 set of indices for which $p_j =1$,
 we have the identity
 $$\int_{\R^M}\prod_{j\in S}^M f_j(\ba_j\cdot x){\rm d}^Mx
=(\det(A_SA_S^t))^{1/2}\prod_{j\in S}^M\int_\R f_j(y){\rm d}y  \ .\Eq(yo1bb)$$
which gives us $D(\br) =  (\det(A_SA_S^t))^{1/2}$. 
This is finite, and since $D(\br)$ is convex and finite at the vertices of $K_A$, it is finite throughout $K_A$.
 \eop

\bigskip
\chap {5: Determination of the optimizers}5
\medskip

A partial solution  to the problem of determining all maximizers, when they exist, for
\eqv(bl1)  was obtained in [\rcite{BL}] where it is proved
that in the case $M =2$ and $N=3$, which gives the classical Young's inequality,  
the only non negative maximizers of the ratio in \eqv(bl1) are certain specific Gaussian functions. 
The method of proof extends to more general cases involving $M+1$ functions in $\R^M$, but not to values of $N > M+1$.  
  
Under the additional assumption that there exists an optimizer to the Gaussian variational problem
\eqv(bl), a full determination of the non negative optimizers was obtained by Barthe [\rcite{Ba2}]. He conjectures that
when there is no optimizer to the Gaussian problem, there is no optimizer for the general problem \eqv(bl1).
Here we give a proof of Barthe's theorem, and of his conjecture. We also determine the form  of all of 
the complex valued optimizers.

Before we begin, note  a  restriction that  we may impose on $\{\ba_1,\ba_2,\dots, \ba_N\}$ without 
loss of generality: We may assume
that if any one vector is deleted from $\{\ba_1,\ba_2,\dots,\ba_N\}$,
then what remains still  spans $\R^M$.  The point is that
when $\ba_1$ is necessary for the whole set to span, there is a change of coordinates  
under which
$$\int_{\R^M}\prod_{j=1}^N f_j(\ba_j\cdot x){\rm d}^Nx = 
{1\over |\bu_1\cdot\ba_1|} \int_{\R}f_1(z_1){\rm d}z_1
\left(\prod_{j=2}^N\int_{\R^{M-1}}f(\bb_j\cdot \bw){\rm d}^{M-1}w\right)\ ,\Eq(change) $$
for some vectors $\bb_2,\dots,\bb_N$ in $\R^{M-1}$. (The calculation in \eqv(change)
is carried out at the end of the Appendix.)
This reduces the analysis of \eqv(yo1N) to an integral of the same type, but with one factor
and one dimension fewer. It also shows that in this case, we must have $p_1 =1$ to obtain 
a finite constant $D$.  Also it is clear in this case that the optimizers need not be Gaussian, since
 $f_1$ can be any integrable function without affecting the value of the ratio in \eqv(bl1).

We therefore make the following definition:
\medskip
\noindent{\bf Definition} Given a spanning set $S = \{\ba_1,\ba_2,\dots,\ba_N\}$ of vectors 
in $\R^M$, we say that $\ba_j$ is {\it essential} in case $S\backslash\{\ba_j\}$
does not span $\R^M$, and we say that $S$ is {\it properly redundant} in case
no vector in $S$ is essential, and moreover, no two vectors in $S$ are proportional.
\medskip

We can always apply the reduction argument given just above to eliminate any
essential vectors. Notice that if $N = M$, every vector is essential, and in fact, we have the identity
$$\int_{\R^M}\prod_{j=1}^N f_j(\ba_j\cdot x){\rm d}^Nx
=(\det(AA^t))^{1/2}\prod_{j=1}^M\int_\R f_j(y){\rm d}y  \ .\Eq(yo1b)$$
For this reason, we are interested mainly in $N>M$.  

We also see right away from \eqv(yo1b) that if $\br$ is a vertex of $K_A$, so that
$M$ of the $L^p$ indices are $1$, and $N-M$ are $\infty$, then we get maximizers in \eqv(bl1)
if and only if we take each of the $L^\infty$ functions to be constant, and there is no restriction on the $L^1$ functions. 
Hence for $\br$ a vertex of $K_A$, the maximizers are far from unique, and need not be Gaussian.

\medskip
\noindent{\bf Lemma 5.1}\ Let $\{\ba_1,\ba_2,\dots,\ba_N\}$ span $\R^M$.
Let $A = [\ba_1,\ba_2,\dots,\ba_N]$ be the $M\times N$ matrix whose $j$-th column is
$\ba_j$. Let $P$ be the orthogonal projection in $\R^N$ onto the image of $A^t$. Then
$\ba_j$ is essential if and only if $P_{j,j} =1$. 
\medskip

\noindent{\bf Proof:} By definition, $\ba_j$ is essential if and only if there do not exist numbers
$u_1,u_2,\dots,u_N$ such that  
$$u_j\ne 0\qquad{\rm  and}\qquad \sum_{k=1}^Nu_k\ba_k = 0\ .\Eq(yo41)$$
Let $\bu$ be the vector in $\R^N$ whose $k$th
entry is $u_k$. Then $\displaystyle{\sum_{k=1}^Nu_k\ba_k = 0}$ is exactly the
condition for $\bu$ to belong to the kernel of $A$. Hence, $\ba_j$ is essential if and only if
$u_j =0$ for each vector $\bu$ in the kernel of $A$. Let $\be_j$ denote the $j$th
standard basis vector in $\R^N$, so that $u_j = \be_j\cdot \bu$. Then, since the
image of $A^t$ is the orthogonal complement of the kernel of $A$, we have that
$\ba_j$ is essential if and only if $\be_j$ belongs to the image of $A^t$. Clearly,
this is the case if and only if $P_{j,j} =1$. \eop
\medskip

\medskip
\noindent{\bf Theorem  5.2 (Barthe)} {\it 
Let $\{\ba_1,\ba_2,\dots,\ba_N\}$ be any properly redundant spanning set, and let 
$\br \in K_A^\circ$. Then  
the solution $S$ of the Euler--Lagrange equations \eqv(bl2) is unique up to a constant multiple.
Moreover,  non-negative functions $f_1 , \dots, f_N$ satisfy 
$$
\int_{\R^M}\prod_{j=1}^N f_j(\ba_j\cdot x){\rm d}^Mx
=   D(\br) \prod_{j=1}\|f_j\|_{p_j} \Eq(yo1c)$$
if and only if there is a number $c>0$ and a vector $\bb\in {\rm Img}(A^t)$ so that for each $j$,
$f_j(y)$ is a multiple of
$$\exp\left(-{c\over 2}\left(s_j^2(y - b_j)^2\right)\right)\ .\Eq(yo68)$$}
\medskip

\noindent{\bf Proof:} Recall the proof of Theorem 3.1. Fix any $t>0$, and let
$h_j(y)$ denote $\ln(f_j(t,y))$. Note that each $f_j(t,y)$ is smooth and strictly positive, so that each $h_j$ is smooth. 
Then by \eqv(ber111) we must have 
$$\int_{\R^M}\left[ \sum_{i, j}
 h'(\ba_i\cdot x)Q_{i,j} h'(\ba_j\cdot x)\right]F(x){\rm d}^M x = 0\ ,\Eq(eqeq)$$
 where  $F(x) = \prod_{j=1}^Nf_j(\ba_j\cdot x)$ and $Q$ is given by \eqv(qfor):
 $Q = S^{-1}(I-P)S^{-1}$ where $P$ is the orthogonal projection onto the image of $(AS)^t = SA^t$. 
\smallskip
 
\noindent{\it Step 1: (Each $h_j$ is a quadratic polynomial)} \ 
 Since $F$ is strictly positive, it follows from \eqv(eqeq)  that
 the vector 
  $$\bv(x) = \left[\matrix{
h_1'(\ba_1\cdot x)\cr
h_2'(\ba_2\cdot x)\cr
\vdots\cr
h_N'(\ba_N\cdot x)\cr}\right]\ ,\Eq(yo32)$$
is such that $S^{-1}\bv(x)$  lies in  ${\rm Img}(SA^t)$ for every $x$. 
This means that
if $\bu$ is any vector in the kernel of $AS$, then $\bu\cdot S^{-1}\bv(x) =0$ for all $x$. 

For any  vector $\bu$ in the kernel of $AS$, we define the function $\phi(x)$
by
$$\phi(x) = \bu\cdot S^{-1}\bv(x)  = \sum_{j=1}^Nu_j s_j^{-1}h'_j(\ba_j\cdot x)\ .$$
Since $\phi$ vanishes identically, 
$\displaystyle{0 = \nabla \phi(x) = \sum_{j=1}^Nu_j s_j^{-1} h''_j(\ba_j\cdot x)\ba_j}$.
This means that for each $x$,
the vector $\bw(x)$ defined by
$$\bw(x) = \left[\matrix{
s_1^{-1}h_1''(\ba_1\cdot x)u_1\cr
s_2^{-1}h_2''(\ba_2\cdot x)u_2\cr
\vdots\cr
s_N^{-1}h_N''(\ba_N\cdot x)u_N\cr}\right]\ ,\Eq(yo33)$$
lies in the kernel of $A$.

We shall first show that $h_1''(\ba_1\cdot x)$ is constant. To do this, write $\ba_1$
as a linear combination of the other vectors $\ba_j$:
$\ba_1 = \sum_{j=2}^N\alpha_j\ba_j$. This is possible since $\ba_1$ is not essential.
There may be many ways of doing this, but we can always
choose one such that a minimal number of the $\alpha$'s are non zero, which we do.
Suppose that
there are exactly $k$ values of $j$, $j_1,j_2,\dots,j_k$ for which $\alpha_j\ne 0$

The vector $\bu = S^{-1}\left[\matrix{1\cr-\alpha_2\cr \vdots\cr -\alpha_N\cr}\right]$
belongs to the kernel of $AS$. We use this vector
$\bu$ in \eqv(yo33) to define $\bw(x)$.  

Now let $\by$ be any vector in $\R^M$ that is orthogonal to $\ba_{j_k}$, but not
orthogonal to $\ba_1$.  Since $w(x)$ lies in the kernel of $A$ for every $x$, so does the
vector we get when we differentiate each component in the $\by$ direction. That is, for each $x$,
$$\left[\matrix{
(\by\cdot\ba_1)s_1^{-1}h_1'''(\ba_1\cdot x)u_1\cr
(\by\cdot\ba_2)s_2^{-1}h_2'''(\ba_2\cdot x)u_2\cr
\vdots\cr
(\by\cdot\ba_2)s_N^{-1}h_N'''(\ba_N\cdot x)u_N\cr}\right]\ $$
lies in the kernel of $A$. 
The $j_k$th component of this vector vanishes identically since $\by\cdot \ba_{{j_k}}= 0$. 
The $N -k$ other entries for which $u_j =0$ also vanish identically.  This means that
for each $x$, the above vector lies in the kernel of $A$, and has no more than $k-1$
non zero entries. By assumption there is no vector in the kernel of $A$ whose first component
is non zero and that has fewer than $k$ non zero entries. Hence the first component must
be zero. Since $\by\cdot\ba_1\ne 0$, and $u_1 \ne 0$, this means
$h_1'''(\ba_1\cdot x) = 0$, and proves that $h_1''$ is constant. 

The argument may now be repeated for each $j$, and we learn at this point that
each $h_j$ is a quadratic function.  

\smallskip
 
\noindent{\it Step 2:  (Determination of  $h''_j$)} \ 
Let $c_j$ denote the value of $h''_j$.  Then, from \eqv(yo33), for any vector $\bu$ in the kernel of $AS$,
the vector whose $j$th entry is $s_j^{-1}c_ju_j$ belongs to ${\rm Ker}(A)$. 
Since $\bu$ is in the kernel of $AS$ if and only if $S\bu$ is in the kernel of $A$, 
we see that $s_j^{-1}c_j$ must be a constant multiple of $s_j$. In other words, for some constant $c$, we have
$$s_1^{-2}h_1''(\ba_1\cdot x) = s_2^{-2}h_2''(\ba_2\cdot x) = \cdots = s_N^{-2}h_N''(\ba_N\cdot x) = -c\Eq(yo31)$$
This of course means that for each $j$, there are constants $a_j$ and $b_j$
so that
$$h_j(y) = -{c\over 2}s_j^2\left(y - b_j\right)^2 + a_j\ ,\Eq(poy3)$$
which means that
$$f_j(y) = \exp\left(-{c\over 2}s_j^2\left(y-b_j\right)^2 + a_j\right)\ .\Eq(yo34)$$
Evidently, $c>0$. 

\smallskip
 
\noindent{\it Step 3:  (Determination of  $b_j$)} \ 
Let $\bb$ denote the vector $\bb$ in $\R^N$ whose $j$th component is $b_j$.
{}From \eqv(poy3) and the definition of $\bv(x)$, we see that $\bv(0) = -cS^{2}\bb$. 
We have seen that $S^{-1}\bv(0)$ lies in ${\rm Img}(SA^t)$, and so $\bb$ lies in 
 ${\rm Img}(A^t)$.

The constant $b_j$ is the mean of the probability density $f_j^{p_j}(y)/\|f_j\|_{p_j}^{p_j}$, and the mean does not
change under the evolution considered here, which commutes with 
translations.  Therefore,
we see that $f_j(t,y)$ can have the form specified in \eqv(yo34) if and only if it has this form at
$t=0$. That is, there is equality in \eqv(yo1c)  if and only if 
there is a positive constant $c$, and a vector $\bb$ in
the image of $A^t$ so that each $f_j$ has the form specified in 
\eqv(yo34) with $b_j$
being the $j$th component of $\bb$, and the $a_j$ are arbitrary.
\eop

 \medskip
\noindent{\bf Corrolary 5.3} {\it 
Let $\{\ba_1,\ba_2,\dots,\ba_N\}$ be any properly redundant spanning set. Then the function 
$\phi(t_1,t_2,\dots,t_N) = Tr\left(\ln\left(Ae^TA^t\right)\right)$ is strictly convex, except along the
lines obtained by adding a number $c$ to each $t_j$. }

\medskip
\noindent{\bf Proof:} Were this not the case, we would have two solutions $S$
of the Euler-Lagrange that would not be constant multiples of one another. \eop

\medskip

The {\it strict}  convexity  was proved by Brascamp and Lieb under the stronger hypothesis that every subset 
of $M$ vectors chosen from  $\{\ba_1,\ba_2,\dots,\ba_N\}$ is linearly independent.

Concerning maximizers for $\br$ on the boundary of $K_A$, we have already dealt  with the
vertices -- these have plenty of non--Gaussian maximizers, and in the strict sense considered here do not
have any Gaussian optimizers:  If $p_j =1$, then  $f_j$ may be any non negative $L^1$ function, and and so may be 
taken to be Gaussian, while if $p_j = \infty$, then $f_j$ must be constant, and therefore not Gaussian.
One could consider constants as degenerate Gaussians, though this would not be entirely consistent with the
terminology we have been using in reference to the Gaussian optimization problem. Alternately, one can stipulate that
$p_j < \infty$ for all $j$. Indeed, if $p_j = \infty$, then the corresponding factors involving $f_j$ 
can be deleted top and bottom in \eqv(bl1) without affecting the constant.   
We may then prove a conjecture of Barthe [\rcite{Ba2}]:
\medskip
\medskip
\noindent {\bf Theorem 5.4} {\it 
Let $\{\ba_1,\ba_2,\dots,\ba_N\}$ be any properly redundant spanning set, and let 
$\br\in K_A\backslash K_A^\circ$ be such that $p_j < \infty$ for all $j$.  Then there may be
no optimizers for the Brascamp Lieb inequality, but when there are optimizers, there are Gaussian optimizers.
Moreover, there is a constructive procedure for deciding whether or not optimizers exist in any particular case.}
\medskip
\noindent{\bf Proof:}  We again apply  the factorization formula \eqv(factor) from Lemma 4.1. 
As in the proof of Theorem 4.3,  let $S$ be a critical set of least cardinality. Such a set exists since  
$\br\in K_A\backslash K_A^\circ$.
As shown in the proof of Theorem 4.3, $D(\br) = D_SD_{S^c}$ where $D_S$ and $D_{S^c}$ are defined in
\eqv(ber201) and \eqv(ber202) respectively.  Since $S$ was a critical subset of least cardinality, there are
 no critical subsets for this problem. Hence there are Gaussian optimizers for the variational problem that determined $D_S$.   
 
 Next, suppose that there are no critical sets in the variational problem that determines $D_{S^c}$.  
 Then this problem has only Gaussian optimizers, unique up to a common scaling and certain translations. However, examining
 \eqv(factor) we see that for 
 $$\int_{\R^r}\prod_{j\in S} f_j(\bb_j\cdot y){\rm d}^r y \left(
\int_{\R^{M-r}}\prod_{j\in S^c} f_j(\bb_j\cdot y + \bc_j\cdot z){\rm d}^{M-r}z\right)\Eq(fac43)$$
to equal $D_SD_{S^c}\prod_{j=1}^N\|f\|_{p_j}$, it is necessary that the translations in the integral on the
right be among those permitted by  Theorem 5.2.  There is a simple criterion for this: Let $A_{S^c}$ be
the matrix obtained by deleting from $A$ the $j$th collumn whenever $j\in S$. Let $U$ and $V$ be the partial 
isometries used in the proof of Lemma 4.1, so that if we put $B = A^t_{S^c}U$ and $C = A^t_{S^c}V$, the columns of $B$ 
(resp. $C$)
are the vectors $b_j$ (resp. $c_j$) in the second integral.

When optimizers exist, it must be the case that for  each $y$ the translation in the second integral is 
one permitted by Theorem 5.2. Clearly this is the case if and only if ${\rm Img}(B) \subset {\rm Img}(C)$. 
Conversely, if this is the case, all of the translation are admissible, and using Gaussian optimizers for $D_{S^c}$
in \eqv(fac43), we will have this integral equal to   $D_SD_{S^c}\prod_{j=1}^N\|f\|_{p_j}$.

The general case is handled in very much the same way:   If there are  critical sets in the variational 
problem that determines $D_{S^c}$, ``peel these off''  repeatedly until one gets a problem 
with no critical subsets, and hence
Gaussian optimizers. Now one works ones way back up, checking the compatibility condition 
${\rm Img}(B) \subset {\rm Img}(C)$  each step of the way.  If this is ever violated, there are no optimizers.
Otherwise, we obtain a set of Gaussian optimizers. \eop
\medskip

One might further hope that the Gaussian functions in Theorem 5.2 are also the only optimizers of
Young's inequality in the wider class of complex valued functions. However, this
is not the case. The reason is that there exist in general functions $\phi_j(y)$ with
$$
e^{i \sum_{j=1}^N \phi_j(\vec a_j \cdot \vec x)} = 1 \ , \Eq(phases)
$$
and thus if $f_1, \dots, f_N$ is any set of non-negative optimizers, then
$e^{\phi_1}f_1, \dots,e^{\phi_N} f_N$ is a set of complex optimizers. 
Here are some examples. Any three vectors in $\R^2$ are linearly dependent, i.e., there
is a relation
$\sum_{j=1}^3 \alpha_j \ba_j =0$. hence with $\phi_j(y) = \alpha_j y$ \eqv(phases) holds.
With four vectors there are more possibilities. E.g., pick $\ba_1= \e_1, \ba_2=\be_2$,
$\ba_3=(\be_1+\be_2)/\sqrt 2$ and  $\ba_4=(\be_1-\be_2)/\sqrt 2$. then the function
$\phi_1(y)=\phi_2(y)=-y^2$, $ \phi_3(y)=\phi_4(y)=y^2$ again satsify \eqv(phases).
Hence there are non-trivial complex valued optimizers. \footnote{*}{\eightpoint 
The possibility of complex optimizers of this type  
for Young's inequality was pointed out to Brascamp and 
Lieb by J. Fournier; see a note added in proof at the end of their  paper.}

In general, let $\{f_1,f_2,\dots,f_N\}$ be any set of optimizers. 
Define functions $z_j$ by
$$z_j(y) = f_j(y)/|f_j(y)|$$
where $f_j(y) \ne 0$, and $z_j(y) =1$ otherwise.
These functions take values in the unit circle in the complex plane. 
In order to have equality in the generalized
Young inequality, it is necessary that 
$$
\prod_{j=1}^N z_j (\vec a_j \cdot \vec x) = 1 \Eq(phases2)$$
almost everywhere.

\medskip
\noindent{\bf Theorem 5.5} {\it Let $\{\ba_1,\dots,\ba_N\}$ be any set of vectors spanning 
$\R^M$ such that no two vectors are multiples of one another. Let $z_j$, $j =1,\dots, N$, 
be any $j$ measurable functions from $\R$ to 
the unit circle in the complex plane such that \eqv(phases2) holds almost everywhere.
Then for each $j$,
$$z_j(y) = e^{i\phi_j(y)}$$
where $\phi_j$  is a polynomial of degree  at most $N-M$. }
\medskip

We first prove a lemma:
\medskip
\noindent{\bf Lemma 5.6} {\it Let $z$ be a function from $\R$ to the unit circle 
in the complex plane,
and let $n$ be any positive integer.  Suppose that $z$ has the following property:
$${z(x+y)\over z(x)}  = e^{i\psi(x,y)}$$
where $\psi(x,y)$ is a polynomial of degree $n-1$ in $x$ with coefficients that a 
measurable functions of $y$.   Then $z(x) = e^{i\phi(x)}$ where $\phi$ is a 
polynomial of degree $n$.}
\medskip

\noindent{\bf Proof:} 
Before beginning, notice that the modulus of $z$ is constant, and non zero. Hence $z$ is never zero.

Consider first the case $n=1$.  Writing $w(y) = e^{\psi(y)}$, 
since there is no $x$ dependence in this case, 
$$z(x+y) = z(x)w(y)\ .$$
Now let $\rho$ be any smooth compactly supported function on $\R$. Then
$$\int z(x+y)\rho(y){\rm d}y = z(x)\int w(y)\rho(y){\rm d}y\ .$$
Since the restriction of $w$ to any interval is a non--zero function  
in $L^2$ on that interval, and since smooth, compactly supported functions 
are dense in this $L^2$ space, we can choose
$\rho$ so that $\int w(y)\rho(y){\rm d}y  = c\ne 0$.  Then we have
$$z(x) = {1\over c}\int z(x+y)\rho(y){\rm d}y\ .$$
This shows that $z$ is smooth.  In particular, once we chose a branch of the 
logarithm for $z(0)$,
there is just one way to choose the logarithm of $z(x)$ so that it is 
continuous, and then of course it is smooth. Hence there is a smooth 
real function $\phi$ so that
$z(x) = e^{i\phi(x)}$, and
$$\phi(x+y) = \phi(x) + \psi(y)\ .$$
Evidently, $\psi$ is also smooth.  Applying $\partial^2/\partial x\partial y$ 
to both sides, we learn
that $\phi''$ vanishes identically, and so $\phi$ is a polynomial of first degree.

Now suppose that $n\ge 2$. Here the argument is similar, but requires one more step. 
We first write
$$z(x+y) = z(x)e^{i\psi(x,y)}\ .$$
Pick any $x_0$, and choose a smooth, compactly supported function $\rho$ as above so that
for this $x_0$,
$$\int  e^{i\psi(x_0,y)} \rho(y){\rm d}y \ne 0\ .$$

Now, no matter how large the coefficients of the polynomial  
$\psi(x,y)$ may be at certain $y$ in the
support of $\rho$, the function
$$x \mapsto   \int  e^{i\psi(x,y)} \rho(y){\rm d}y$$
is continuous in $x$ by the Dominated Convergence Theorem.

We conclude that $c(x) = \int  e^{i\psi(x_0,y)} \rho(y){\rm d}y$ 
is continuous and non--zero on a neighborhood of $x_0$. Hence
$$z(x) = {1\over c(x)}\int z(x+y)\rho(y){\rm d}y$$
is continuous on a neighborhood of $x_0$. Since $x_0$ is arbitrary, $z$ is continuous.

It now follows that   $e^{i\psi(x,y)}$ is continuous in both $x$ and $y$. 
Therefore, the coefficients
are uniformly  bounded functions of $y$ in any compact interval. 
This means that all of the
partial derivatives in $x$ of   $e^{i\psi(x,y)}$ are integrable and continuous, 
and so the function
$c(x) =   \int  e^{i\psi(x,y)} \rho(y){\rm d}y$ that we defined above is not only
continuous, it  is actually smooth in $x$. It now follows that $z$ is smooth, and as before
we write $z = e^{i\phi}$, and have
$$\phi(x+y) = \phi(x) + \psi(x,y)\ .$$
Taking $n$ derivatives in $x$, and using the hypothesis that $ \psi(x,y)$ 
has degree $n-1$
in $x$, we see that the $n$th derivative of $\phi$ is constant. 
Hence $\phi$ is a polynomial of degree $n$. \eop
\medskip

It is of course well known that if $\phi$ and $\psi$ are two measurable functions on 
$\R$ such that
$$\phi(x+y) = \phi(x) + \psi(y)$$
then both $\phi$ and $\psi$ are first degree polynomials. 
Lemma 5.6 generalizes this in several respects. 
It seems likely that it may be known, but we cannot  find any reference for it.
\medskip

\noindent{\bf Proof of Theorem 5.5:} We can easily eliminate any essential  vectors from $\{\ba_1,\dots,\ba_N\}$:
If $\ba_j$ is essential, it is clear that $z_j$ is constant. Hence we may assume that
$\{\ba_1,\dots,\ba_N\}$ is properly spanning.   

It suffices by symmetry to show that $z_1$ has the specified form.
Choose a basis for $\R^M$ from   $\{\ba_1,\dots,\ba_N\}$ that contains $\ba_1$. 
After renumbering, we may assume this is   $\{\ba_1,\dots,\ba_M\}$.  Let $\bb_1$ be unit vector
that is orthogonal to the span of  $\{\ba_2,\dots,\ba_M\}$, and scaled so that  $\bb_1\cdot \ba_1 =1$.

Now, for any $y_1$ in $\R$, translate the identity \eqv(phases2) by replacing $x$ with $x+y_1\bb_1$. 
Since $\bb_1$ is orthogonal 
to $\ba_j$ for $2\le j \le M$, the corresponding factors are unaffected by translation, and hence
$$z_1((\vec a_1 \cdot \vec x)+ y_1)(\prod_{j=M+1}^N z_j (\vec a_j \cdot (\vec x+ y_1\bb_1)) = 
z_1 (\vec a_1 \cdot \vec x)\prod_{j=M+1}^N z_j (\vec a_j \cdot \vec x)\ .$$

Let $T_y$ be the operator
$$T_y(z)(x) = {z(x+y)\over z(x)}\ .$$ Then defining 
$$z_1^{(1)}(t;y_1) = {T_{y_1}z_1\over z_1}(t)\ ,$$ and defining
$w_j = z_j (\vec a_j \cdot (\vec x+ y_1\bb_1))/z_j (\vec a_j \cdot \vec x)$ for $j \ge M+1$, 
$$z_1^{(1)}(\vec a_1\cdot x;y_1)\prod_{j=M+1}^N w_j (\vec a_j \cdot \vec x) =1\ .\Eq(phases7p)$$
This is of the same form as \eqv(phases2), but with fewer functions. 

Next, choose $\bb_2$ so that $\bb_2\cdot\ba_{M+1} = 0$ (if it wasn't the case already
that $\bb_1\cdot \ba_{M+1} =0$), but $\bb_2\cdot \ba_1 =1$.  Making the same sort of translation in \eqv(phases7p),
but this time by $y_2\bb_2$,
we eliminate the second factor by dividing through, so that the first factor becomes
$$z_1^{(2)}(t;y_1,y_2) = {T_{y_2}z_1^{(1)}\over z_1^{(1)}}(t;y_1)\ .$$

Proceeding in this way, we eventually learn that for some $k < M-N$,
$${T_{y_{k+1}}z_1^{(k)}\over z_1^{(k)}}(t;y_1,\dots,y_k)$$
is independent of $t$. 

By Lemma 5.6, it follows that $z_1^{(k)}(t;y_1,\dots,y_k) = e^{i\phi(t;y_1,\dots,y_k)}$
where $\phi(t;y_1,\dots,y_k)$ is a first degree polynomial in $t$ with coefficients that are measurable in
$y_1,\dots,y_k$. But by definition,
$$z_1^{(k)}(t;y_1,y_2) = {T_{y_k}z_1^{(k-1)}\over z_1^{(k-1)}}(t;y_1,\dots,y_{k-1})\ .$$
Applying Lemma 5.6 again, we learn the form of $z_1^{(k-1)}$. Proceeding in this way, we learn the form of
$z_1$. \eop
\medskip

Once one knows that the possible phase functions are polynomials of limited degree, it is a problem in 
linear algebra to determine them explicitly for any particular set of vectors
$\{\ba_1,\dots,\ba_N\}$.

\bigskip
\chap {6: The best best constant }6
\medskip

Let   $ \{\ba_1,\ba_2,\dots,\ba_N\}$ be a properly redundant spanning set
 in $\R^M$, and let  $P$ denote the orthogonal projection onto the image of $A^t$. 
Notice that
$Tr(P) = {\rm rank}(A^t) = M$.
Also, since $P$ is an orthogonal projection, each diagonal entry $P_{j,j}$ satisfies
$0\le P_{j,j} \le 1$. Furthermore, since no column of $A$ is zero, we actually have
$0 < P_{j,j} $
for each $j$. 

Indeed, since ${\rm rank}(A) = M$, $AA^t$ is positive definite, and
$P = A^t(AA^t)^{-1}A$. Therefore, 
$$P_{j,j} = \be_j\cdot  A^t(AA^t)^{-1}A\be_j = \ba_j\cdot  (AA^t)^{-1} \ba_j > 0\ .$$
Hence, if we define $p_j^{\phantom{.}\circ}$ by
$${1\over p_j^{\phantom{.}\circ}} = P_{j,j} = \ba_j\cdot  (AA^t)^{-1} \ba_j > 0\ ,$$
we have that whenever $ \{\ba_1,\ba_2,\dots,\ba_N\}$
is properly redundant, $1 <  p_j^{\phantom{.}\circ} < \infty$ for each $j$, and also 
$\sum_{j=1}^N(1/p_j) =  Tr(P) = M$, so that \eqv(yo22) is satisfied.
Morover, the Euler--Lagrange equation \eqv(bl2) is then satisfied with $S = I$. 

\medskip
\noindent{\bf Definition} Let $ \{\ba_1,\ba_2,\dots,\ba_N\}$  be a properly redundant
spanning set of vectors 
in $\R^M$, and let  $A = [\ba_1,\ba_2,\dots,\ba_N]$. Let $P$ be the orthogonal
projection in $\R^N$ onto the image of $A^t$. For $j= 1,2,\dots,N$,
define $p^{\phantom{.}\circ}_j = 1/P_{j,j}$. Then $\br^{\phantom{.}\circ} = 
\{p^{\phantom{.}\circ}_1,p^{\phantom{.}\circ}_2,\dots,p^{\phantom{.}\circ}_N\}$ is the 
{\it canonical set of $L^p$ indices corresponding to} $ \{\ba_1,\ba_2,\dots,\ba_N\}$.   The terminology will be justified
by Theorem 6.1 below.

\medskip

Since  for $\br = {\br}^{\phantom{.}\circ}$,  the Euler--Lagrange equations 
\eqv(bl2) are satisfied with $S = I$, it follows from \eqv(ber45) that 
$$D(\br^{\phantom{.}\circ}) = 
\prod_{j=1}^N\left( P_{j,j}\right) ^{-P_{j,j}/2}{\rm det}(AA^t)^{-1/2}\ .\Eq(ber62)$$

Notice that while computing $D(\br)$ for given $L^p$ indices is a  nonlinear optimization problem,
calculating the  $D(\br_0)$ is a simple matter of linear algebra. This is significant since it 
turns out that given the vectors
$\{\ba_1,\ba_2,\dots,\ba_N\}$,
$D(\br^{\phantom{.}\circ})$ is the ``best best constant'' in the generalized Young's inequality
$$\int_{\R^M}\prod_{j=1}^N f_j(\ba\cdot x){\rm d}^Nx \le 
D\prod_{j=1}^N\|f_j\|_{p_j}$$
This  justifes the terminology ``canonical $L^p$ indices'':

\medskip
\noindent{\bf Theorem 6.1:} {\it For any properly redundant spanning set,
$D(\br^{\phantom{.}\circ}) < D(\br)$
for all $\br \ne \br^{\phantom{.}\circ}$. }

\medskip
\noindent{\bf Proof:} Let $\phi_A$ be the function defined by
\eqv(yo63), and $\phi_A^*$ its Legendre transform.
It was shown by Brascamp and Lieb that    $\phi_A$ is convex. Since $\phi_A$ is smooth 
as well as convex, $\phi_A^*$ is strictly convex.  From \eqv(yo65) and \eqv(yo65b), we have
$$2\ln(\tilde D(\br )) 
=\phi_A^*\left({1\over p_1},{1\over p_2},\dots, {1\over p_N}\right)\Eq(yo65c)$$

By the definition of $\phi_A$, and the Euler--Lagrange equation
\eqv(bl2), if $1/p^{\phantom{.}\circ}_j$ is the $j$th canonical $L^p$ index,
$$\left({1\over p^{\phantom{.}\circ}_1},{1\over p^{\phantom{.}\circ}_2},\dots, 
{1\over p^{\phantom{.}\circ}_N}\right) = \nabla \phi_A(0)\ .\Eq(nab1)$$
But since the gradients of Legendre transforms are inverse to one another,
$$\nabla\phi_A^*(\nabla\phi_A(0)) = 0\ .$$
This proves that the vector on the left in \eqv(nab1) is a critical point of $\phi_A^*$. Since
$\phi_A^*$ is strictly convex, it is the unique minimizer. \eop
\medskip

We also note that formula \eqv(yo65c) displays $D(\br)$ as a log--convex function of $\br$, This can be used to 
produce {\it arbitrarily sharp} upper bounds on $D(\br)$ for a given set of $L^p$ indices: Using Newton's method
or some other means of generating explicit approximate solutions of the Euler--Lagrange 
equations \eqv(bl2), generate several approximate solutions. For each, compute the ``best best constant'' 
for each $\{s_1\ba_1,s_2\ba_2,\dots,s_N\ba_N\}$.  If $\br$ can be written 
as a convex combination of the corresponding vectors of cannionical inverse $L^p$ indices,
then $D(\br)$ can be bounded above by a convex combination of the corresponding ``best best constants''.

Special cases of the canonical $L^p$ indices have arisen in applications of the Brascamp Lieb inequality. 
A beautiful application to convex geometry by Keith Ball [\rcite{Ball}] concerened a situation in which $N$ unit vectors
$\vec u_1,\dots, \vec u_N$ satisify 
$$\sum_{j=1}^N c_j\vec u_j\vec u_j^t = I_{M\times M}\ .\Eq(ball)$$
where the $c_j$ are positive numbers. Clearly, $\sum_{j=1}^N c_j = M$. Let $\vec a_j = \sqrt{c_j}\vec u_j$. Then
\eqv(ball) becomes $AA^t = I_{M\times M}$. It follows that the orthogonal projection onto the image of $A^t$
is simply $A^tA$, and the $j$th diagonal entry is $c_j$. Hence taking $p_j = 1/c_j$ gives the canonical
$L^p$ indices in this case. These were the $L^p$ indices used by Ball in his application. 

Since for the canonical $L^p$ indices, 
the Euler--Lagrange equation \eqv(bl2) is then satisfied with $S = I$,
the heat flow interpolation argument of Section 3 gives an especially simple proof of the 
inequality in this case. For this reason, the
method of proof developed here works very simply in Keith Ball's context; 
see [\rcite{BaCE}] for more information.

\medskip

\noindent{\bf Example:} Consider the five vectors
$$\ba_1 = \left[\matrix{\phantom{-}1\cr  -1\cr \phantom{-}0 \cr}\right] \qquad
\ba_2 = \left[\matrix{\phantom{-}0\cr  \phantom{-}1\cr -1 \cr}\right] \qquad
\ba_3 = \left[\matrix{-1\cr   \phantom{-}0 \cr \phantom{-}1}\right] \qquad
\ba_4 = \left[\matrix{1\cr  0\cr 0 \cr}\right] \qquad
\ba_5 = \left[\matrix{0\cr 1\cr 0 \cr}\right] \ .$$
It is easily seen that this is a properly redundant spanning set. Notice that the first three
vectors all lie in the plane $x_1+x_2+x_3 =0$. As long as $0 <1/ p_j < 1$ for each $j$,
$$\sum_{j\in S}{1\over p_j} < |S|\ .\Eq(yo87)$$ 
Therefore, as long as $r(S) \ge \min\{|S|,M\}$ and \eqv(yo22) is satisfied, and there are no
supercritical sets.  The only set $S$ with
$r(S) <  \min\{|S|,M\}$
is $S = \{1,2,3\}$. Therefore, as long as
${\displaystyle {1\over p_1} + {1\over p_1} + {1\over p_1}  < 2}$,
together with  the scaling condition  \eqv(yo22) and $0 < p_j < 1$ for each $j$ and are all satisfied, 
 $\br$ belongs to $K_A^\circ $, and $K_A$ is the closure of
 the points obtained in this way.  An easy computation shows that the canonical indices for this example are $p_1 = 2$,
and $p_2 = p_3 = p_4 = p_5 =8/5$.  By Theorem 4.4, $K_A$ is has  $9$ vertices, and is their convex hull.
\medskip

\bigskip
\chap {Appendix}7 
\medskip

In this section we exhibit trial functions that show the optimality of Theorem 1
and Theorem 2, and describe the change of variables leading to \eqv(change).  

First we show that the inequality in Theorem 1 cannot hold with {\it any} constant  if the index 
$p$ of the $L^p$ norms on the right side is less than $2$.  
  For any given $0 < \alpha < 1$, define  $f(v)$ be defined by
$$f(v) = |v|^{-\alpha} + (1-v^2)^{-\alpha(N-1)/2}\ .\Eq(trial)$$
Then $\displaystyle{\int_{[-1,1]} f^p {\rm d}\nu_N < \infty}$ as long as $p\alpha < 1$,
as one easily sees from \eqv(marg).

On the other hand, discarding one term in each factor, 
$$\prod_{j=1}^N f(v_j) \ge \left(\prod_{j=1}^{N-1} |v_j|^{-\alpha}\right)(1 - v_N^2)^{-\alpha(N-1)/2}\ .$$
We can parameterize the upper and lower hemispheres of $S^{N_1}$ using the coordinates $(v_1,\dots,v_{N-1})$
The intergal over $S^{N-1}$ is then easily converted into an integral over the unit ball in $\R^{N-1}$. 
Doing this in radial coordinates, we have, since $|v_n| = \sqrt{1 - r^2}$ in these coordinates,
$$\int_{S^{N-1}}\left(\prod_{j=1}^N f(v_j)\right){\rm d}\mu \ge C\int_0^1 r^{-2\alpha(N-1)}{r^{N-2}\over \sqrt{1- r^2}}
{\rm d}r$$
where $C$ is a positive constant resulting from the angular integration. This integral diverges unless $\alpha < 1/2$.

The conclusion is that for all $N$ and all $p<2$, there is a positive function $f$ so that
${\displaystyle\int_{S^{N-1}}f^p(v_1){\rm d}\mu < \infty}$
while
${\displaystyle\int_{S^{N-1}}\left(\prod_{j=1}^N f(v_j)\right){\rm d}\mu = \infty}$.

Next we turn to the entropy inequality in Theorem 2. Consider a spherical cap on the sphere $S^{N-1}$ centered at 
the point $v_1=1, v_2=0, \dots, v_N=0$
of radius $\varepsilon$ denote its characteristic function by $\chi_\varepsilon$. 
Define 
$$F  = H\chi_\varepsilon \qquad{\rm with}\qquad  
H= \left(\int_{S^{N-1}} \chi_\varepsilon {\rm d} \mu\right)^{-1}\ .$$
Clearly, $H$ is of order
$\varepsilon^{-(N-1)}$, and hence for $\varepsilon$ small, $S(F)$ is of order
$$
-\log(H) \Eq(largeentropy)
$$
which is of order $(N-1) \log(\varepsilon)$. Since the function is invariant under all rotation that
fix the $v_1$ axis, we get that the entropy of the marginal is also given by \eqv(largeentropy).
Moreover, the $j$-th marginal can be thought of as averaging the function $H\chi_\varepsilon$
over all roations that keep the axis $v_j$ fixed. the resulting function is essentially a 
multiple of a characteristic
function of a band of width $2 \varepsilon$ that is centered at the equator perpendicular to the $v_j$ axis. Call this 
function $f_j:=\psi_\varepsilon$. Since
the integral of this function must be equals to one the height of this function must be
$h= (\int_{S^{N-1}} \psi_\varepsilon {\rm d} \mu)^{-1}$,
and is of order $1/\varepsilon$. Hence its entropy is of order
$\log(\varepsilon)$. Thus the sum of the entropies of the marginals is given, in leading order,  by
$2(N-1) \log(\varepsilon)$
which is twice the entropy of the function $F$. This shows that the constant $2$ in the
entropy inequality is sharp.

Finally, the coordinate change leading to \eqv(change)
may be described as follows:
Suppose that $\ba_1$ is not  in the span of 
$\{\ba_2,\dots,\ba_N\}$.  Let $\{\bu_1,\bu_2,\dots,\bu_M\}$ be an orthonormal basis of 
$\R^M$ so that  $\{\bu_2,\dots,\bu_M\}$ has the same span as 
$\{\ba_2,\dots,\ba_N\}$. Let $R$ be the matrix given by $R = [\ba_1,\bu_2,\dots,\bu_M]$.
(That is, the first column of $R$ is $\ba_1$, the second column is $\bu_2$, and so forth).
Then $R$ is invertible, and we can define new coordinates $z$ by $z = R^tx$. With
this definition, $z_1 = \ba_1\cdot x$.  Moreover, for $j\ge 2$,
$$\ba_j\cdot x = \ba_j\cdot (R^t)^{-1}z = (R^{-1}\ba_j)\cdot z\ .$$
Since $R^{-1}\ba_j$ is the coordinate vector of $\ba_j$ with respect to the basis
$\{\ba_1,\bu_2,\dots,\bu_M\}$, $(R^{-1}\ba_j)_1 =0$ for $j\ge 2$.
Therefore, defining $\bw$ in $\R^{M-1}$ by $w_j = z_{j+1}$, there are uniquely determined vectors
$\bb_j$ in $\R^{M-1}$ so that
$(R^{-1}\ba_j)\cdot z = \bb_j\cdot\bw$. Since 
$${\rm d}^Mx = {1\over |\bu_1\cdot\ba_1|}{\rm d}^Mz = 
{1\over |\bu_1\cdot\ba_1|}{\rm d}z_1{\rm d}^{M-1}w\ ,$$
we have \eqv(change).

\bigskip

\noindent{\bf Bibliography}

\bigskip

\item{[\rtag{Ball}]} {K. Ball:}
``Volume Ratios and a reverse isoperimetric inequality'' 
Jour. London Math. Soc., {\bf 44} no. 2 351--359 (1991)

\vskip 2mm

\item{[\rtag{Ba}]} {F. Barthe:}
``Optimal Young's inequality and its converse, a simple proof'' 
Geom. Func. Analysis., {\bf 80} 234--242 (1998)

\vskip 2mm

\item{[\rtag{Ba2}]} {F. Barthe:}
``On a reverse form of the Brascamp--Lieb inequality'' 
Invent. Math., {\bf 134}  no. 2,  235--361 (1998)

\vskip 2mm

\item{[\rtag{BaCE}]} {F. Barthe and D. Cordero--Erausquin:}
``Inverse Brascamp--Lieb inequalities along the heat equation'' 
to appear in {\it Geometric Aspects  of Functional Analysis, 2002--2003},
eds. V. Milman and G. Schechtman, Lecture Notes in Mathematics 1850, Springer Verlag,
Berlin 2004

\vskip 2mm

\item{[\rtag{BL}]} {H. J. Brascamp and  E. H. Lieb:}
``Best constants in Young's inequality, its converse, and its generalization to 
more than three functions'' Advances in Math., {\bf 20} 151-173 (1976)

\vskip 2mm

\item{[\rtag{CCL1}]} {E. A. Carlen, M. C. Carvalho and  M. Loss:}
Many body aspects of approach to equilibrium,
in Journes Equations aux derivees partielles, Nantes, 5-9 juin 2000.

\vskip 2mm

\item{[\rtag{CCL2}]} E. A. Carlen, M. C. Carvalho, M. Loss, ``Determination of the spectral
gap for Kac's master equation and related stochastic evolution'',
Acta Mathematica {\bf 191}, 1-54 (2003). arXiv math-ph/0109003

\vskip 2mm

\item{[\rtag{L}]} E. H. Lieb, ``Gaussian kernels have only Gaussian maximizers'',
Invent Math {\bf 102}, 179--208 (1990). 

\vskip 2mm

\end